\documentclass[format=acmsmall,nonacm]{acmart}

\pdfoutput=1 % For arxiv
\usepackage[utf8]{inputenc}
\usepackage[ruled]{algorithm2e}
\renewcommand{\algorithmcfname}{ALGORITHM}
\SetAlFnt{\small}
\SetAlCapFnt{\small}
\SetAlCapNameFnt{\small}
\SetAlCapHSkip{0pt}
\IncMargin{-\parindent}
\setcopyright{none}
\authorsaddresses{} % None please
\acmDOI{}
\acmArticle{}
\usepackage{amsmath} 		% math
\usepackage{myshortcuts}	% personal shortcuts
\usepackage{jlcode}
\usepackage{bm}

\acmVolume{0X}
\acmNumber{0X}
\acmArticle{0X}
\acmYear{2018}
\acmMonth{0X}

\begin{document}

% Page heads
%\markboth{G. Zhou et al.}{A Multifrequency MAC Specially Designed for WSN Applications}

% Title portion
\title{NEP-PACK: A Julia package for nonlinear eigenproblems\\
Release v0.2}
\author{Elias Jarlebring}
\affiliation{\institution{KTH Royal Institute of Technology}}
\author{Max Bennedich}
\affiliation{\institution{KTH Royal Institute of Technology}}
\author{Giampaolo Mele}
\affiliation{\institution{KTH Royal Institute of Technology}}
\author{Emil Ringh}
\affiliation{\institution{KTH Royal Institute of Technology}}
\author{Parikshit Upadhyaya}
\affiliation{\institution{KTH Royal Institute of Technology}}

%\author{GANG ZHOU
%\affil{College of William and Mary}
%YAFENG WU
%\affil{University of Virginia}
%TING YAN
%\affil{Eaton Innovation Center}
%TIAN HE
%\affil{University of Minnesota}
%CHENGDU HUANG
%\affil{Google}
%JOHN A. STANKOVIC
%\affil{University of Virginia}
%TAREK F. ABDELZAHER
%\affil{University of Illinois at Urbana-Champaign}}
% NOTE! Affiliations placed here should be for the institution where the
%       BULK of the research was done. If the author has gone to a new
%       institution, before publication, the (above) affiliation should NOT be changed.
%       The authors 'current' address may be given in the "Author's addresses:" block (below).
%       So for example, Mr. Abdelzaher, the bulk of the research was done at UIUC, and he is
%       currently affiliated with NASA.

\begin{abstract}
  We present NEP-PACK a novel open-source library for
  the solution of nonlinear eigenvalue problems (NEPs). The package
  provides a framework to represent NEPs, as well as efficient
  implementations
  of many state-of-the-art algorithms. The package
  makes full use of the efficiency of Julia, yet maintains
  usability,  and integrates well with other software packages.
  The package is designed to be easy to use for
  application researchers as well as algorithm developers.
  Particular attention is paid to algorithm neutrality, in
  order to make performance comparisons between algorithms
  easier. This paper describes the main functionality
  of NEP-PACK, as well as design decisions and
  theory needed for the design.
\end{abstract}

%\category{C.2.2}{Computer-Communication Networks}{Network Protocols}

%\terms{Design, Algorithms, Performance}

%\keywords{Wireless sensor networks, media access control,
%multi-channel, radio interference, time synchronization}

%\acmformat{Gang Zhou, Yafeng Wu, Ting Yan, Tian He, Chengdu Huang, John A. Stankovic,
%and Tarek F. Abdelzaher, 2010. A multifrequency MAC specially
%designed for  wireless sensor network applications.}
% At a minimum you need to supply the author names, year and a title.
% IMPORTANT:
% Full first names whenever they are known, surname last, followed by a period.
% In the case of two authors, 'and' is placed between them.
% In the case of three or more authors, the serial comma is used, that is, all author names
% except the last one but including the penultimate author's name are followed by a comma,
% and then 'and' is placed before the final author's name.
% If only first and middle initials are known, then each initial
% is followed by a period and they are separated by a space.
% The remaining information (journal title, volume, article number, date, etc.) is 'auto-generated'.

%\begin{bottomstuff}
%\end{bottomstuff}

\maketitle

\renewcommand{\shortauthors}{Jarlebring et al.}

\medskip\medskip
\section{Introduction}

This package concerns nonlinear eigenvalue problems defined as the problem
of determining the singular points of a matrix, i.e., find $(\lambda,v)$ such that
\begin{equation}  \label{eq:nep}
  M(\lambda)v=0
\end{equation}
where $v\neq 0$ and $M:\CC\rightarrow \CC^{n\times m}$ is a holomorphic (or meromorphic
with only a few poles).

Nonlinear problems which are not the linear or generalized eigenvalue problem,
occur in many situations. Some of the most common situations are
\begin{itemize}
\item higher order differential equations (references in \cite{Tisseur:2001:QUADRATIC}) leading to matrix polynomials
\item systems and control for time-delay systems, leading to exponential nonlinearities \cite{Michiels:2007:STABILITYBOOK} \cite{Jarlebring:2008:THESIS}
\item quantum physics (quantum dots) \cite{Betcke:2007:THESIS,Betcke:2007:QUANTUM}) leading to rational nonlinear functions
\item fluid mechanics (scaling exponent in turbulent flow) \cite{Ooi:2017:NEP} leading to exponential nonlinearities
\item fluid-solid interations \cite{Voss2005,Voss:2003:MAXMIN} leading to rational terms,
\item boundary element method applied to resoncance problems \cite{Steinbach:2007:BEM,Steinbach:2009:BEM,unger2013convergence}, see also \cite{Effenberger2012a} and the softare package \cite{Smigaj:2015:BEM},
\item absorbing boundary conditions (frequency dependent) leading
  to square root nonlinearities \cite{Tausch2000,Tausch2002,Jarlebring:2017:TIAR} or
  Bessel functions  \cite{ARAUJO:2018:Helmholtz}, e.g., in
  fiber optics design \cite{Kaufman2006},
\item chatter in machine tool milling \cite{Bueler:2007:ERROR,Insberger:2002:MATHIEU,Insperger:2002:SEMIDISC,Jarlebring:2017:BROYDEN,Rott:2010:ITERATIVE},
\item periodic structures, e.g., in crystals \cite{SCHMIDT:2009:BAND,Engstrm2010,Effenberger2012,Fliss:2017:DTN}
\end{itemize}

In most of these applications there is need for performance,
and robustness. Our package is implemented in the Julia programming language
\cite{Bezanson2017}, in order to obtain efficiency and still have
access to high-level functionality. A milestone for computing
in the Julia language was carried out within the Celeste project,
which qualifies as petascale computation \cite{celeste2018}.

%Nonlinear eigenvalue problem algorithms yields an interesting
%show-case for Julia, since many algorithms involve
%several levels of loops (low-level performance)
%and still require access to other
%efficient matrix functions (high-level functionality).
%

%\begin{center}
%Generality \& easy to use   \;\;\;$\Leftrightarrow$\;\;\; Efficency \& performance
%\end{center}
%
%\begin{center}
%Beginner   \;\;\;$\Leftrightarrow$\;\;\; Advanced user
%\end{center}
%
% We need a algorithm neutral way to represent NEPs.

% Elevator pitch: As easy and efficient use as \verb#full(), sparse, and speye()# for matrices

%Software manuscripts we may want to use as style-inspiration:
%\begin{itemize}
%\item DDE-BIFTOOL: \cite{Engelborghs:2002:DDEBIFTOOL}
%%\item BEM++ \cite{Smigaj:2015:BEM}
%\item The SIAM rev. Julia publication: \cite{Bezanson2017}
%\end{itemize}
%
The numerical treatment of this problem has received attention in a large
number of works, see summary papers
such as \cite{Ruhe:1973:NLEVP}, \cite{Mehrmann:2004:NLEVP}
\cite{Voss2012} and \cite{Guttel2017} as
well as software packages \cite{Betcke2010}
and \cite{Roman:2018:SLEPC,Hernandez:2003:SSL,Hernandez:2005:SSF}.

As we shall further describe in Section~\ref{sect:representation},
many applications and algorithms are based on a sum of products representation
of the $M$ matrix. We will provide considerable functionality and
efficiency for problems that can be expressed as
\begin{equation} \label{eq:SPMF_NEP}
  M(\lambda)=A_1f_1(\lambda)+\cdots+A_mf_m(\lambda).
\end{equation}
In theory, any NEP can be expressed as \eqref{eq:SPMF_NEP}, if one sets $m=n^2$.
However, most algorithms based on \eqref{eq:SPMF_NEP} also assume
that $m$ is not too large, and become less attractive
due to an increase in computation time unless  $m$ is relatively small.
Our framework is efficient for structures as
\eqref{eq:SPMF_NEP}.

Our softare is designed to not be based on \eqref{eq:SPMF_NEP}
but rather on certain interface functions, which define a NEP. This
allows to represent NEPs, where \eqref{eq:SPMF_NEP} is not efficient. These framework
interfaces are described in Section~\ref{sect:interfaces}.

All of the NEP-algorithms are carefully documented in terms of references,
and in order to encourage users to give credit to the original
algorithm researchers.

\section{Basic usage}
\begin{figure}[h]
  \begin{center}
%\pgfplotsset{grid style={dotted,gray}} % change grid style (minor/major grid style)
\begin{tikzpicture}
  %\begin{tikzpicture}[domain=0:4]% if we want to plot without axis
  %\begin{tikzpicture}[y=.2cm, x=.7cm,font=\sffamily] % set size
  \begin{axis}
    [width=8.0cm, height=6cm,% 6cm normally enough space to use as two subfigures
      /pgf/number format/.cd,
      1000 sep={},
      xlabel={Re $\lambda$},ylabel={Im $\lambda$},%
      xmin=-1.8,xmax=1,%
      %ymin=1E-8,ymax=1E3,%
      grid,%
      %ticks=none, %
      %axis lines = middle, axis line style={->, line width=0.8pt},
      %ytick={10^{-16},1,...,10},%
      %minor ytick={1,2,...,10},%
      legend entries={ ~$\lambda$},%
      %legend pos=north east,
      %legend style={at={(0.5,-0.4)},anchor=north}, % below plot
      %legend pos=outer north east, % new version of pgf
      %legend style={draw,cells={anchor=west},row sep=0pt},%
      %axis x line*=bottom, % only line at bottom
      %mark repeat={3},% skip drawing the mark every X timme
      %axis equal image,
    ]
    \addplot[only marks, color=blue] coordinates {
( -1.62378  ,  9.3967      )
( -1.34453  ,  9.2985      )
( -1.62378  , -9.3967      )
( -1.34453  , -9.2985      )
( -0.883271 ,  5.32505     )
( -1.296    ,  5.05755     )
( -0.455848 ,  1.68846     )
( -0.097128 ,  2.47396e-13 )
(  0.308866 , -6.54814e-13 )
( -0.455848 , -1.68846     )
( -1.296    , -5.05755     )
( -0.883271 , -5.32505     )};
    %\addplot[red, solid, only marks, mark options={solid,scale=0.5}] table[x={x}, y={y}]{data.dat}; % data stored as columns
  \end{axis}
\end{tikzpicture}
    \caption{The eigenvalues of the neuron example.
      \label{fig:neuron}
    }
  \end{center}
\end{figure}
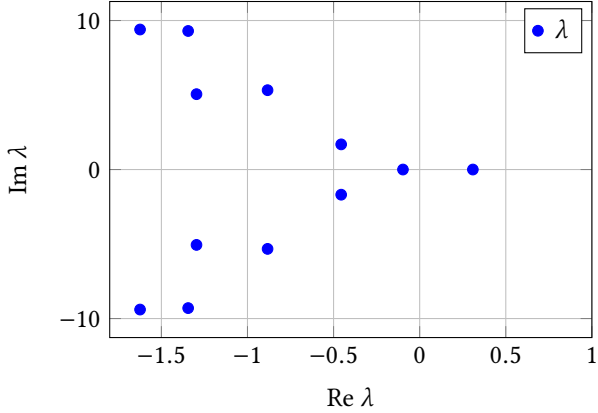

NEP-PACK is a registered package in the Julia central package repository, which
makes it possible to install the package with very little effort
\begin{lstlisting}
julia> ]
(v1.0) pkg> add NonlinearEigenproblems
julia> using NonlinearEigenproblems
\end{lstlisting}
Nonlinear eigenvalue problems are represented as objects of
the type NEP, which can be created in a number of different ways.
We have a gallery of problems available, which can be accessed through
the \verb#nep_gallery# command (further described in Section~\ref{sec:benchmark}).
\begin{lstlisting}
julia> nep=nep_gallery("neuron0");
\end{lstlisting}
This creates a NEP object which is used to model a neuron.
This gallery problem stems from \cite{Shayer:2000:STABILITY}
which is also available as a model problem in DDE-BIFTOOL
\cite{Engelborghs:2002:DDEBIFTOOL,Engelborghs:2001:DDEBIFTOOL}. The problem describes the stability
 the delay differential equation
%\begin{subequations}
\begin{eqnarray}
\dot{x}_1(t)&=&-\kappa x_1(t)+\beta\tanh(x_1(t-\tau_3))+a_1\tanh(x_2(t-\tau_2))\\
\dot{x}_2(t)&=&-\kappa x_2(t)+\beta\tanh(x_2(t-\tau_3))+a_2\tanh(x_1(t-\tau_1)).
\end{eqnarray}
%\end{subequations}
%In particular, the stability of the zero solution is given by the solution to the NEP with
In particular, the stability of the zero solution is characterized by the eigenvalues with the largest real part of the following NEP
\[
M(\lambda):=-\lambda I+A_0+A_1e^{-\tau_1\lambda}+A_2e^{-\tau_2\lambda}+A_3e^{-\tau_3\lambda},
\]
%which is commonly called a delay eigenvalue problem (DEP).
which belongs to the class of NEP commonly called a delay eigenvalue problems (DEP), see Section~\ref{sec:types}.
The typeof command reveals that the problem is represented as a \verb#DEP#:
\begin{lstlisting}
julia> typeof(nep)
DEP{Float64,Array{Complex{Float64},2}}
\end{lstlisting}
As an illustration we solve this problem
with our implementation of the NLEIGS method \cite{Guttel:2014:NLEIGS}
\begin{lstlisting}
julia> Σ=[-3.0-10im,-3+10im,1+10im,1-10im]; # Region of interest
julia> (λ,V)=nleigs(nep,Σ)
julia> using Plots; plotly();
julia> plot(λ,linewidth=0, markershape=:xcross, label="λ")
\end{lstlisting}
The same results are obtained with the infinite Arnoldi method \cite{Jarlebring:2010:DELAYARNOLDI}
\begin{lstlisting}
julia> (λ,V)=iar_chebyshev(nep,σ=-2) # keyword argument σ = target point
\end{lstlisting}
The result is given in Figure~\ref{fig:neuron}.
%\begin{figure}[h]
%  \begin{center}
%    %\subfigure[mytext]{\scalebox{0.65}{\includegraphics{file}}}
%    \includegraphics[width=5cm]{neuronexample.png}
%    \caption{The eigenvalues of the neuron example.
%      \label{fig:neuron}
%    }
%  \end{center}
%\end{figure}
%%julia> Ω=[-1.0-1im,-1+1im,1+1im,1-1im];
%%julia> nleigs(ComplexF64,nep,Ω)
%

%\begin{filecontents}{data.dat} % move to preamble
%  x    y
%  2000  3E-7
%  2005  3E-4
%\end{filecontents}

\section{The NEP representation}\label{sect:representation}
\subsection{Accessing the data in a NEP}\label{sect:interfaces}
The vast amount of numerical methods and the vast amount of applications
that have been  developed and formulated in the literature have
typically been  expressed in the way that
is the most natural for that specific method/application.
Different methods require data from the NEP in different ways. Different
applications are expressed in different ways.
Our package is designed for algorithms and problems
which can be (efficiently) expressed with
the following quantities, which are selected
to match the need from as many applications and algorithms
as possible.

\begin{itemize}
\item \verb#compute_Mder(lambda,k)# computes the $k$th derivative
  of $M$ and evaluates it in $\lambda$
  \begin{equation}  \label{eq:mder}
 M_k=M^{(k)}(\lambda)
  \end{equation}
\item \verb#compute_Mlincomb(lambda,V)# computes the linear
  combination
  \begin{equation}  \label{eq:mlincomb}
    \sum_{i=1}^kM^{(i-1)}(\lambda)v_i.
  \end{equation}
 given the evaluation point $\lambda$ and a matrix $V\in\CC^{n\times k}$.
\item \verb#compute_MM(S,V)# computes the expression
  \begin{equation}  \label{eq:MMdef}
\MM(S,V)=\frac{1}{2\pi i}\int_\Gamma M(\xi)V(\xi I-S)^{-1}\,d\xi
  \end{equation}
  for given matrices  $S\in\CC^{p\times p}$  and $V\in\CC^{n\times p}$,
  and the contour $\Gamma$ includes the eigenvalues of $S$. The form \eqref{eq:MMdef} is more commonly expressed in terms
of matrix functions
\[
   \MM(S,V)=A_1Vf_1(S)+\cdots+A_mVf_m(S),
   \]
  if $A_1,\ldots,$ and $f_1,\ldots$ are given as in \eqref{eq:SPMF_NEP}.
\end{itemize}

\noindent Note that these
compute functions are mathematically equivalent, i.e., there are explicit
procedures to compute one quantity from any other quantity. We specify
further relations in Section~\ref{sec:equivalences}. Although they are
mathematically equivalent, they are computationally very different and
the transformations are not necessarily very efficient.
For instance, in many applications (particular types of matrix-free situations) the matrix
may not be directly available, but only available as a subroutine.
In this case the user can specify \eqref{eq:mlincomb},
whereas computing a matrix as in \eqref{eq:mder} may not
be computationally feasable.

Routines for some additional secondary quantities are also available:
\begin{itemize}
\item \verb#lin_solve# solves a linear system associated with $M(\lambda)$. Further specifications can be made by inheriting from the type \verb#LinSolver#.
\item \verb#compute_rf# computes the Rayleigh functional for the NEP, i.e., solves
  the nonlinear (scalar) equation
  \begin{equation}  \label{eq:rf}
    y^HM(\lambda)x=0.
  \end{equation}
\item Several algorithms require the solution to a projected problem
  \begin{equation}  \label{eq:proj}
   Y^HM(\lambda)Xz=0
  \end{equation}
  which is again a NEP, see Section~\ref{sec:projection}. The fact that projected
  NEPs are also of the type NEP, allows us to apply any of our methods as an inner
  solver.
\item \verb#errmeasure# is a keyword argument accepted by most functions. The function handle should accept two arguments \verb#errmeasure(lambda,v)# and computes an error estimate based on $\lambda$ and $v$. The \verb#default_errmeasure# computes the relative residual norm
 \[
 \frac{\| M(\lambda)v\|}{\|v\|}.
 \]
 The construction allows the user to specify in which way the error should be measured.
 Hence other error measurements, such as the backward error presented in \cite{Higham:2008:BACKWARD}, can be implemented by the user.
\end{itemize}
In what follows we describe in what way state-of-the-art algorithms can be implemented
with these compute functions.
\begin{itemize}
\item A large class of methods can be derived from Newton's metod, e.g., \cite{Neumaier:1985:RESINV,Unger:1950:NICHTLINEARE,Schreiber:2008:PHD,Spence:2005:PHOTONIC,Ruhe:1973:NLEVP,Anselone:1968:NEWTON,Jarlebring:2018:DISGUISED}. See also the summary of papers in, e.g., \cite{Szyld:2015:LOCAL,Guttel2017,Szyld2013,Szyld:2011:EFFICIENT,Szyld2014}. These algorithms have in common that one needs
  to compute for some vectors $u_1$ and $u_2$
  \[
    M(\lambda)u_1+M'(\lambda)u_2
    \]
    where in some settings $u_1$ is zero. This is clearly possible with \eqref{eq:mlincomb}. They
    also require the solution to one (or many) linear systems, which can be computed with the
    \verb#lin_solve# functionality.
  \item The standard application of contour integral methods require
    the solution to many linear systems (for different evaluation points) and
    only require a matrix vector product, i.e., \eqref{eq:mlincomb}. Methods in these papers
    can be classified in this way: \cite{Asakura:2010:NUMINTPEP,Asakura:2009:NUMERICAL,Beyn:2011:INTEGRAL} and the FEAST software \cite{Polizzi:2009:FEAST}, and acceleration techniques \cite{Xiao:2017:SOLVING}.
  \item Jacobi-Davidson methods \cite{Betcke:2004:JD,Sleijpen1996,Voss:2007:JD,Effenberger2013} are projection methods and hence require \eqref{eq:proj}. For the deflation technique presented in \cite{Effenberger2013}, the quantity
  \[U(\mu)=\frac{1}{2\pi}\int_\Gamma M(\xi)V(\xi I - S)^{-1}(\xi-\mu)^{-1} \,d\xi\]
  and its derivatives are required in order for the deflated problem to implement the interfaces desired by the solvers of \eqref{eq:proj}. These derivatives can be computed via the interface \eqref{eq:mlincomb} and the observations that $U^{(k)}(\mu) = (-M^{(k)}(\mu)V + k U^{(k-1)}(\mu))(S-\mu I)^{-1}$, and $U(\mu) = (\MM(S,V) - \MM(\mu I,V))(S-\mu I)^{-1}$. Specifically, due to structure the action of  $\MM(\mu I,V)$ can be computed using \eqref{eq:mlincomb}, and for an invariant pair $\MM(S,V)=0$.
\item Other projection methods such as the nonlinear Arnoldi method \cite{Voss:2004:ARNOLDI} and variation \cite{Jarlebring:2005:RATIONAL} require the projection \eqref{eq:proj}. Also the block preconditioned PCG for NEPs in \cite{Xue:2018:BHP} is a projection method. When solving large-scale problems, the storage of the projection subspace becomes costly as the number of iterations increase. In such cases, we need restarting techniques as discussed in \cite{BETCKE:2017:RESTART}, which can
  be directly implemented with manipulations of the projected NEP.
 \item Infinite Arnoldi type methods such as \cite{Jarlebring:2010:DELAYARNOLDI,Jarlebring:2010:TAYLORARNOLDI,Jarlebring:2012:INFARNOLDI} and two-sided version \cite{Gaaf:2017:INFBILANCZOS} require \eqref{eq:mlincomb}. The restarted versions \cite{Mele2018,Jarlebring:2017:SCHUR} additionally require \eqref{eq:MMdef}.
 \item Rayleigh functional methods, e.g., \cite{Neumaier:1985:RESINV,Schreiber:2008:PHD,Schwetlick:2012:NONLINEAR} depend on solution methods for the Rayleigh quotient.
 \item Certain algorithms are based on directly working with the block residual
   \eqref{eq:MMdef}, e.g., block Newton \cite{Kressner:2009:BLOCKNEWTON,Effenberger2012,Effenberger2012a,Beyn2010,Beyn2011,Effenberger2013a}
 \item Methods based on the QR-method, e.g., \cite{Garrett2016,Kublanovskaya:1970:APPROACH} require \eqref{eq:mder}
 \item Rational Krylov methods such as, e.g., \cite{VanBeeumen:2013:RATIONAL,Guttel:2014:NLEIGS,VanBeeumen:2015:COMPACTRATIONAL} requires access to linear solvers. Moreover, for a general NEP either \eqref{eq:mder} or \eqref{eq:mlincomb} is used, and in the special case of an SPMF, see Section~\ref{sec:types}, it uses matrix function evaluations $f_0(S),\dots,f_m(S)$.
\end{itemize}
\subsection{The common types} \label{sec:types}
We provide the user with efficient implementations
of the compute functions in the previous section,
for many common types:
\begin{itemize}
\item \verb#PEP#: $M(\lambda)=A_1+\lambda A_2+\cdots \lambda^{m-1}A_{m}$
\item \verb#DEP#: $M(\lambda)=-\lambda I+A_0+\sum_{i=1}^me^{-\tau_i\lambda}A_i$
\item \verb#SumNEP#: $M(\lambda)=A(\lambda)+B(\lambda)$ where $A$ and $B$ are also \verb#NEP#s
\item \verb#SPMF_NEP#: See below.
\item \verb#LowRankNEP#: An SPMF where the matrices are represented as low-rank factorizations
\end{itemize}
The most general of the above is the \verb#SPMF_NEP# which
represents the sum of products of matrices and functions
\eqref{eq:SPMF_NEP}.
 The functions $f_i$, $i=1,\ldots,m$ have to be defined in scalar sense, as well as in a matrix function sense.
 The implementation of the SPMF is designed to be efficient when $m\ll n$.
 The example below solves the NEP
 \[
  M(\lambda)=\lambda A+e^{\lambda}B+\left(1+\sqrt{\lambda}\right)C
  \]
  with the block Newton method \cite{Kressner:2009:BLOCKNEWTON}:
 \begin{lstlisting}
julia> using LinearAlgebra
julia> A=ones(5,5); B=ones(5,5)+I; C=reverse(B,dims=1)
julia> f1= S-> S;
julia> f2= S-> exp(S)
julia> f3= S-> one(S)+sqrt(S);
julia> nep=SPMF_NEP([A,B,C],[f1,f2,f3]);
julia> blocknewton(nep,S=[1 0; 0 1.0], X=[1 0; 0 1; zeros(3,2)],displaylevel=1)
Iteration 1: Error: 2.112578e+01
Iteration 2: Error: 3.705499e+00
Iteration 3: Error: 2.361554e+00
Iteration 4: Error: 2.955760e-01
Iteration 5: Error: 3.543752e-03
Iteration 6: Error: 4.080188e-07
Iteration 7: Error: 2.885914e-15
(Complex{Float64}[0.557832+0.0im -3.03756e-16+0.0im; -7.18644e-16+0.0im
  0.557832+0.0im], Complex{Float64}[-0.617521+0.0im -0.00206428+0.0im;
  0.00206428-0.0im -0.617521+0.0im; ... ; -0.00206428+0.0im 0.617521-0.0im;
  0.617521-0.0im 0.00206428-0.0im])
 \end{lstlisting}
\subsection{Equivalence of interfaces}\label{sec:equivalences}
As an example how \eqref{eq:mlincomb} can be computed from \eqref{eq:MMdef}:
$\MM(S,V) e_1$ is equal to \eqref{eq:mlincomb} where $S \in \RR^{k \times k}$ is the bidiagonal matrix with $\lambda$ in the main diagonal, $S_{i+1,i}=i$ and $V=[v_1, \dots, v_k]$.
This equivalence follows by expressing $M(\lambda)$ in \verb#SPMF_NEP# format \eqref{eq:SPMF_NEP} in \eqref{eq:mlincomb} and by using \cite[Definition 1.2]{Higham:2006:FUNCTIONSOFMATRICES} with a proper rescaling.

\section{Problem transformations}
The abstraction of the NEPs to be essentially specified by well-defined
compute functions (Section~\ref{sect:representation})
leads to the advantage that problems
can be transformed leading by defining new compute functions.

We have implemented a number of ways to transform the problem
\begin{itemize}
\item One can shift and scale the problem, i.e., define a new NEP
  \begin{align} \label{eq:M_shift_scale}
    \tilde{M}(\lambda)=M(\alpha\lambda+\sigma)
  \end{align}
 This functionality is available in the \verb#shift_and_scale# function.
\item One can carry out a Möbius transformation of the problem, i.e.,
  define a new NEP
  \[
 \tilde{M}(\lambda)=M((a\lambda+b)/(c\lambda+d))
 \]
 This functionality is available in the \verb#mobius_transformation# function.
\item One can deflate eigenvalue (or invariant pairs) from a NEP
  as specified, e.g., in \cite{Effenberger2013a}. This is provided by
  the function \verb#effenberger_deflation#
\end{itemize}
Although the above functions provide convenient features for a user, they
may not always lead to extremely efficient algorithms, since a transformed problem
may have some computational overhead. Therefore, certain functionality is also provided
at an algorithm level, e.g., shifting and scaling is available in
the infinite Arnoldi methods.

The deflation can be used to compute one pair at a time and avoid reconvergence, e.g.,
as follows:
\begin{lstlisting}
julia> nep=nep_gallery("dep0");
julia> (s,v)=newton(nep);
julia> n=size(nep,1);
julia> S0=reshape([s],1,1);
julia> V0=reshape(v,n,1);
julia> dnep=effenberger_deflation(nep,S0,V0)
julia> (s2,v2)=augnewton(dnep);  # this converges to different eigval
julia> minimum(svdvals(compute_Mder(nep,s2)))
9.323003321058995e-17
\end{lstlisting}
\section{NEP-Solver algorithm implementations}
We have implemented several algorithms as well as extensions.
\subsection{Newton-type methods}
Several flavors of Newton's method are available.
Armijo rule steplength combined with deflation increases reliability
of these methods considerably.

\begin{itemize}
\item \verb#augnewton#: Augmented Newton \cite{Unger:1950:NICHTLINEARE}
\item \verb#resinv#: Residual inverse iteration \cite{Neumaier:1985:RESINV}
\item \verb#blocknewton#: Block Newton method \cite{Kressner:2009:BLOCKNEWTON}
\item \verb#quasinewton#: Quasi-Newton method \cite{Jarlebring:2018:DISGUISED}
\item \verb#implicitdet#: Implicit determinant \cite{Spence:2005:PHOTONIC}
\item \verb#newtonqr#: Newton QR approach \cite{Kublanovskaya:1970:APPROACH}
\item \verb#mslp#: Method of successive linear problems \cite{Ruhe:1973:NLEVP}
\item \verb#sgiter#: Safe-guarded iteration \cite{Voss:2004:ARNOLDI}
\item \verb#rfi#: Rayleigh functional iteration \cite{Schreiber:2008:PHD}
\item \verb#broyden#: Broydens method \cite{Jarlebring:2017:BROYDEN}
\end{itemize}
\subsection{Krylov-based methods}
\begin{itemize}
  \item \verb#nlar# Nonlinear Arnoldi method \cite{Voss:2004:ARNOLDI}
  \item \verb#nleigs# NLEIGS \cite{Guttel:2014:NLEIGS}
  \item \verb#iar# (and variants) infinite Arnoldi method \cite{Jarlebring:2012:INFARNOLDI,Jarlebring:2017:TIAR}
\end{itemize}
\subsection{Projection methods}\label{sec:projection}
\begin{itemize}
\item \verb#nlar# Nonlinear Arnoldi method \cite{Voss:2004:ARNOLDI}
\item \verb#jd# Jacobi-Davidson method \cite{Betcke:2004:JD} and \cite{Effenberger2013}
\end{itemize}
\subsection{Contour integral methods}
\begin{itemize}
  \item \verb#beyn_contour# \cite{Beyn2011}
\end{itemize}
%\subsection{Inner-outer iteration construction}
%
%JD + NLArnoldi + IAR (optionally) can use inner-outer iteration.
%Design carefully so we can use same system structure for all.
%
%Features (generality, extendibility, efficiency, easy to use)
%\begin{itemize}
%\item Efficiency from multiple dispatch
%\item Generality and extendability by the possibility to define new types,
%\item Easy to use - since we have a number of inner solvers directly available and can be changed with keyword argument.
%\end{itemize}
%
\section{Benchmark problems}\label{sec:benchmark}
We have made a number of
benchmark problems available via the \verb#nep_gallery# command, e.g.,
a standardized delay eigenvalue problem can be
loaded with
\begin{lstlisting}
julia> nep=nep_gallery("dep0");
\end{lstlisting}
Several large-scale problems, such as the model
of the waveguide in \cite{Jarlebring:2017:TIAR,Ringh:2018:SYLVPRECOND} are available.

The library of  Berlin-Manchester benchmark problem in the MATLAB NLEVP package \cite{Betcke:2013:NLEVP},
can be accessed in NEP-PACK in two ways. A subset of the problems from that collection
have been converted to native NEP-PACK format, e.g., the
``gun'' problem can be loaded with the commands
\begin{lstlisting}
julia> nep=nep_gallery("nlevp_native_gun");
\end{lstlisting}
Several implementation techniques had to be adapted to Julia in
order to become efficient, e.g., the Bessel function nonlinearity
in the ``fiber'' benchmark as described in Section~\ref{sec:fiber}.
The NLEVP problems can also be accessed by using the  Julia packages
which can communicate with a MATLAB process running in the background. We have provided wrappers
such that the problems can be loaded with the command:
\begin{lstlisting}
julia> using GalleryNLEVP
julia> nep=nep_gallery(NLEVP_NEP,"fiber")
julia> quasinewton(nep,λ=1e-6)
(7.139494342432901e-7 + 5.123670257712833e-18im, Complex{Float64}[-97388.3
  -10508.4im,-2.75452e5-29721.8im, -5.06025e5-54601.1im, -7.79049e5-84060.9im,
  -1.0887e6-1.17473e5im, -1.43106e6-1.54414e5im, -1.80322e6-194571.0im,
  -2.20294e6-237702.0im, -2.62842e6-2.83611e5im, -3.07814e6-3.32137e5im
  ...  -3.67741e7-3.968e6im, -3.6741e7-3.96443e6im, -3.6708e7-3.96087e6im,
  -3.6675e7-3.95731e6im, -3.66421e7-3.95376e6im, -3.66092e7-3.9502e6im,
  -3.65763e7-3.94666e6im, -3.65434e7-3.94311e6im, -3.65106e7-3.93957e6im,
  -3.64778e7-3.93603e6im])
\end{lstlisting}
Note that the wrapper is completely transparent such that \verb#quasinewton# makes
a call to the NLEVP library (available in a MATLAB process which runs in the background) every time
it accesses the NEP. Due to the communication overhead, it is generally preferred
to use the native methods for larger problems, due to the overhead generated
by the communication between Julia and MATLAB.
\section{Performance comparison}

\subsection{NLEIGS Julia implementation}\label{sec:comparison_nleigs}
We want to provide empirical support for the performance
of our package, and the Julia language.
In order to do so, we used the MATLAB NLEIGS implementation\footnote{NLEIGS version 0.5 available for download at \url{http://twr.cs.kuleuven.be/research/software/nleps/nleigs.html}.}
described in \cite{Guttel:2014:NLEIGS}. For illustration
purposes we converted the MATLAB code to Julia and the
NEP-PACK procedures to access data,
such that it can be considered a good candidate to assess the performance of
Julia vs MATLAB.

%We have compared NEP-PACK version 0.2.1 with

We used the same two large scale problems as in the above paper and
included as benchmarks in the MATLAB implementation;
the "gun" problem, and the "particle in a canyon" problem, and we ran
the same six experiments. See \cite{Guttel:2014:NLEIGS} for
full details. The experiments were run on a MacBook Pro,
with a 2.9 GHz Intel i7-6920HQ,
2x4 cores, and 16 GB memory. We used MATLAB v8.4.0 (R2014b) and Julia v1.0.2.
The results are reported in Table~\ref{tbl:comparison}. Our Julia reimplementation
is faster and consumes less memory. The general explanation is the way
Julia handles data structures, which improves the possibility
to carry out Just-In-Time compilation. More precisely,
we observed that the inner loops (often consisting of orthgonalization)
were considerably faster, also handling of sparse matrices differed
considerably in performance.

Each MATLAB experiment was carried out 20 times,
and the fastest run is reported. For the Julia implementation we used
the Benchmark toolbox, with parameter \verb#seconds=500#. We report the median CPU-time
for the NEP-PACK implementation in Table~\ref{tbl:comparison}.
The memory usage is the amount of memory used at the end of the algorithm,
including cached LU factors. Note that although the implementations should
behave identically, the number of iterations required for convergence
may vary a bit due to different start vectors and tiny rounding errors that
build up over time.

\begin{table}[h]
  \begin{center}
    \begin{tabular}{c|c|c|c|c||c|c|c|c}
      %      & t1 & t2 & t3 \\
      &\multicolumn{4}{|c||}{MATLAB}&
      \multicolumn{4}{|c}{Julia / NEP-PACK}   \\
      \hline
%      Experiment
      & Iter& Conv. $\lambda$&CPU&Memory& Iter& Conv. $\lambda$&CPU &Memory\\
%|----------|:----:|:-------:|:------:|:----:|:----:|:-------:|:------:|:----:|
   \hline
Gun P&100&17&6.4 s&420 MB&100&17&3.9 s&59 MB\\
Gun R1&100&21&6.9 s&421 MB&100&21&4.0 s&59 MB\\
Gun R2&95&21&20.1 s&413 MB&95&21&12.7 s&51 MB\\
Gun S&70&21&5.2 s&408 MB&71&21&3.7 s&46 MB\\
Particle R2&78&2&16.6 s&213 MB&74&2&7.5 s&73 MB\\
Particle S&141&2&13.0 s&239 MB&134&2&5.9 s&92 MB\\

    \end{tabular}
    \caption{Performance comparison of NLEIGS implementation in NEP-PACK and the original MATLAB implementation\label{tbl:comparison}}
%    \caption{Performance comparison}
%    of the MATLAB implementation
%      of NLEIGS and our Julia converted version.
%      %\label{tbl:comparison}
%    }
  \end{center}
\end{table}
%

%%% Old data
%\begin{table}[h]
%  \begin{center}
%    \begin{tabular}{c|c|c|c|c||c|c|c|c}
%      %      & t1 & t2 & t3 \\
%      &\multicolumn{4}{|c||}{MATLAB}&
%      \multicolumn{4}{|c}{Julia / NEP-PACK}   \\
%      \hline
%%      Experiment
%      & Iter& Conv. $\lambda$&CPU&Memory& Iter& Conv. $\lambda$&CPU &Memory\\
%%|----------|:----:|:-------:|:------:|:----:|:----:|:-------:|:------:|:----:|
%   \hline
%Gun P&100&17&7.1 s&420 MB&100&17&4.3 s&59 MB\\
%Gun R1&100&21&7.7 s&421 MB&100&21&4.3 s&59 MB\\
%Gun R2&95&21&22.7 s&413 MB&95&21&14.1 s&51 MB\\
%Gun S&70&21&6.0 s&408 MB&71&21&3.9 s&46 MB\\
%Particle R2&78&2&19.7 s&213 MB&74&2&8.7 s&73 MB\\
%Particle S&141&2&15.1 s&239 MB&134&2&6.5 s&92 MB\\
%
%    \end{tabular}
%    \caption{\label{tbl:comparison}}
%%    \caption{Performance comparison}
%%    of the MATLAB implementation
%%      of NLEIGS and our Julia converted version.
%%      %\label{tbl:comparison}
%%    }
%  \end{center}
%\end{table}
\medskip\medskip
\subsection{Computation of many derivatives}
In order to show the extendability of our framework,
we now show an unusual NEP with $200$ terms. It can be
created and solved as follows.
\begin{lstlisting}
julia> using Random, BenchmarkTools
julia> Random.seed!(0)
julia> m=200;
julia> fv=Vector{Function}(undef,m);
julia> for i=1:m; fv[i]=(x-> exp(i^(1/6)*x)); end;
julia> fv[1]=x->one(x); fv[2]=x->x;
julia> Av=Vector{SparseMatrixCSC}(undef,m);
julia> n=50;
julia> for i=1:m; Av[i]=sprand(n,n,0.01); end;
julia> nep=SPMF_NEP(Av,fv);
julia> v0=ones(n);
julia> @btime iar(nep,maxit=100,v=v0)
  7.568 s (13079550 allocations: 3.88 GiB)
\end{lstlisting}
Due to the fact that the problem has many exponential
terms, the evaluation of the derivatives required in the
infinite Arnoldi method becomes dominant.
Precomputation of derivates are available through the
\verb#DerSPMF#-type, which essentially precomputes derivatives
in a given point, but otherwise behaves as the
parent NEP. The following code shows the improvement.
\begin{lstlisting}
julia> dnep=DerSPMF(nep,0.0,100);
julia> @btime iar(dnep,maxit=100,v=v0)
  3.365 s (12254494 allocations: 1.07 GiB)
\end{lstlisting}
Note that \verb#DerSPMF# is extending the functionality of
standard NEPs, by allowing a precomputation to take place,
but maintain all other functionality of the original NEP.
The \verb#DerSPMF# is again a \verb#NEP# and
precomputation in several points can be achieved by successive
application of the \verb#DerSPMF#.
%
%
%julia> @btime iar(dnep,maxit=100,v=v0)
%  876.589 ms (2557045 allocations: 391.38 MiB)
%
%julia> nep=nep_gallery("nlevp_native_fiber");
%julia> (A0,A1,A2)=get_Av(nep);
%julia> (f0,f1,f2)=get_fv(nep);
%julia> nA0=norm(A0); nA1=norm(A1); nA2=norm(A2)
%julia> errmeasure=(λ,v) -> norm(compute_Mlincomb(nep,λ,v))/(abs(f0(λ))*nA0+abs(f1(λ))*nA1+abs(f2(λ))*nA2)
%julia> γ=1e-5;
%julia> snep=shift_and_scale(nep,scale=γ);
%julia> serrmeasure=(λ,v) -> errmeasure(λ*γ,v);
%julia> λ,v=iar(snep,σ=0.2,errmeasure=serrmeasure,maxit=100,Neig=5);
%julia> using BenchmarkTools
%julia> v0=ones(size(nep,1));
%julia> @btime λ,v=iar(snep,σ=0.2,errmeasure=serrmeasure,maxit=100,Neig=5,v=v0,displaylevel=1);
%  677.176 ms (55774 allocations: 846.04 MiB)
%julia> dsnep=DerSPMF(snep,0.2,40);
%julia> @btime λ,v=iar(dsnep,σ=0.2,errmeasure=serrmeasure,maxit=100,Neig=5,v=v0);
%
\section{Conclusions}
We have presented and described the current release of the package NEP-PACK. The current
state of the software is ready to be used for many use-cases, e.g.,
comparison of algorithms and development of new algorithms, as we have shown in Section~\ref{sec:comparison_nleigs} that it already outperforms other publicly available implementations
of NEP-solvers. Several implementations currently do not have the full functionality,
e.g., some functions do not return eigenvectors but only eigenvalues, although they are available
in theory.
Further testing of other NEP-types, applications
and algorithms to obtain improvements of efficiency for large-scale problems.
The package has been tested on the HPC-environment at KTH Royal Institute of Technology,
and results will be reported in a later version. The development of this package
is done in a public GIT-HUB repository\footnote{\url{https://github.com/nep-pack/NonlinearEigenproblems.jl}} and has a public users manual\footnote{\url{https://nep-pack.github.io/NonlinearEigenproblems.jl}}, in order to improve
possibilities to interact with users and other developers.

\begin{acks}
The authors wish to thank Antti Koskela (Univ. Helsinki)
for comments, suggestions and discussions regarding an early version of the
package. We also wish to express our thanks to researchers who have
made their software available online, e.g., C. Effenberger and R. Van Beeumen.
\end{acks}

% Bibliography
\bibliographystyle{ACM-Reference-Format}
\bibliography{fulljabref}

%%% -*-BibTeX-*-
%%% Do NOT edit. File created by BibTeX with style
%%% ACM-Reference-Format-Journals [18-Jan-2012].

\begin{thebibliography}{82}

%%% ====================================================================
%%% NOTE TO THE USER: you can override these defaults by providing
%%% customized versions of any of these macros before the \bibliography
%%% command.  Each of them MUST provide its own final punctuation,
%%% except for \shownote{}, \showDOI{}, and \showURL{}.  The latter two
%%% do not use final punctuation, in order to avoid confusing it with
%%% the Web address.
%%%
%%% To suppress output of a particular field, define its macro to expand
%%% to an empty string, or better, \unskip, like this:
%%%
%%% \newcommand{\showDOI}[1]{\unskip}   % LaTeX syntax
%%%
%%% \def \showDOI #1{\unskip}           % plain TeX syntax
%%%
%%% ====================================================================

\ifx \showCODEN    \undefined \def \showCODEN     #1{\unskip}     \fi
\ifx \showDOI      \undefined \def \showDOI       #1{#1}\fi
\ifx \showISBNx    \undefined \def \showISBNx     #1{\unskip}     \fi
\ifx \showISBNxiii \undefined \def \showISBNxiii  #1{\unskip}     \fi
\ifx \showISSN     \undefined \def \showISSN      #1{\unskip}     \fi
\ifx \showLCCN     \undefined \def \showLCCN      #1{\unskip}     \fi
\ifx \shownote     \undefined \def \shownote      #1{#1}          \fi
\ifx \showarticletitle \undefined \def \showarticletitle #1{#1}   \fi
\ifx \showURL      \undefined \def \showURL       {\relax}        \fi
% The following commands are used for tagged output and should be
% invisible to TeX
\providecommand\bibfield[2]{#2}
\providecommand\bibinfo[2]{#2}
\providecommand\natexlab[1]{#1}
\providecommand\showeprint[2][]{arXiv:#2}

\bibitem[\protect\citeauthoryear{Anselone and Rall}{Anselone and Rall}{1968}]%
        {Anselone:1968:NEWTON}
\bibfield{author}{\bibinfo{person}{P. Anselone} {and} \bibinfo{person}{L.
  Rall}.} \bibinfo{year}{1968}\natexlab{}.
\newblock \showarticletitle{The solution of characteristic value-vector
  problems by {N}ewton's method}.
\newblock \bibinfo{journal}{\emph{Numer. Math.}}  \bibinfo{volume}{11}
  (\bibinfo{year}{1968}), \bibinfo{pages}{38--45}.
\newblock


\bibitem[\protect\citeauthoryear{Araujo-Cabarcas, Engström, and
  Jarlebring}{Araujo-Cabarcas et~al\mbox{.}}{2018}]%
        {ARAUJO:2018:Helmholtz}
\bibfield{author}{\bibinfo{person}{Juan~Carlos Araujo-Cabarcas},
  \bibinfo{person}{Christian Engström}, {and} \bibinfo{person}{Elias
  Jarlebring}.} \bibinfo{year}{2018}\natexlab{}.
\newblock \showarticletitle{Efficient resonance computations for {Helmholtz}
  problems based on a {Dirichlet-to-Neumann} map}.
\newblock \bibinfo{journal}{\emph{J. Comput. Appl. Math.}}
  \bibinfo{volume}{330} (\bibinfo{year}{2018}), \bibinfo{pages}{177 -- 192}.
\newblock
\showISSN{0377-0427}
\urldef\tempurl%
\url{https://doi.org/10.1016/j.cam.2017.08.012}
\showDOI{\tempurl}


\bibitem[\protect\citeauthoryear{Asakura, Sakurai, Tadano, Ikegami, and
  Kimura}{Asakura et~al\mbox{.}}{2009}]%
        {Asakura:2009:NUMERICAL}
\bibfield{author}{\bibinfo{person}{J. Asakura}, \bibinfo{person}{T. Sakurai},
  \bibinfo{person}{H. Tadano}, \bibinfo{person}{T. Ikegami}, {and}
  \bibinfo{person}{K. Kimura}.} \bibinfo{year}{2009}\natexlab{}.
\newblock \showarticletitle{A numerical method for nonlinear eigenvalue
  problems using contour integrals}.
\newblock \bibinfo{journal}{\emph{JSIAM Letters}}  \bibinfo{volume}{1}
  (\bibinfo{year}{2009}), \bibinfo{pages}{52--55}.
\newblock


\bibitem[\protect\citeauthoryear{Asakura, Sakurai, Tadano, Ikegami, and
  Kimura}{Asakura et~al\mbox{.}}{2010}]%
        {Asakura:2010:NUMINTPEP}
\bibfield{author}{\bibinfo{person}{J. Asakura}, \bibinfo{person}{T. Sakurai},
  \bibinfo{person}{H. Tadano}, \bibinfo{person}{T. Ikegami}, {and}
  \bibinfo{person}{K. Kimura}.} \bibinfo{year}{2010}\natexlab{}.
\newblock \showarticletitle{A numerical method for polynomial eigenvalue
  problems using contour integral}.
\newblock \bibinfo{journal}{\emph{Japan J. Indust. Appl. Math.}}
  \bibinfo{volume}{27} (\bibinfo{year}{2010}), \bibinfo{pages}{73--90}.
\newblock


\bibitem[\protect\citeauthoryear{Betcke}{Betcke}{2007}]%
        {Betcke:2007:THESIS}
\bibfield{author}{\bibinfo{person}{M. Betcke}.}
  \bibinfo{year}{2007}\natexlab{}.
\newblock \emph{\bibinfo{title}{Iterative projection methods for symmetric
  nonlinear eigenvalue problems with applications}}.
\newblock \bibinfo{thesistype}{Ph.D. Dissertation}. \bibinfo{school}{Technical
  University Hamburg-Harburg}.
\newblock


\bibitem[\protect\citeauthoryear{Betcke and Voss}{Betcke and Voss}{2007}]%
        {Betcke:2007:QUANTUM}
\bibfield{author}{\bibinfo{person}{M. Betcke} {and} \bibinfo{person}{H. Voss}.}
  \bibinfo{year}{2007}\natexlab{}.
\newblock \showarticletitle{Stationary {Schrödinger} equations governing
  electronic states of quantum dots in the presence of spinorbit splitting}.
\newblock \bibinfo{journal}{\emph{Appl. Math.}}  \bibinfo{volume}{52}
  (\bibinfo{year}{2007}), \bibinfo{pages}{267 -- 284}.
\newblock


\bibitem[\protect\citeauthoryear{Betcke and Voss}{Betcke and Voss}{2017}]%
        {BETCKE:2017:RESTART}
\bibfield{author}{\bibinfo{person}{Marta~M. Betcke} {and}
  \bibinfo{person}{Heinrich Voss}.} \bibinfo{year}{2017}\natexlab{}.
\newblock \showarticletitle{Restarting iterative projection methods for
  Hermitian nonlinear eigenvalue problems with minmax property}.
\newblock \bibinfo{journal}{\emph{Numer. Math.}} \bibinfo{volume}{135},
  \bibinfo{number}{2} (\bibinfo{date}{01 Feb} \bibinfo{year}{2017}),
  \bibinfo{pages}{397--430}.
\newblock
\showISSN{0945-3245}
\urldef\tempurl%
\url{https://doi.org/10.1007/s00211-016-0804-3}
\showDOI{\tempurl}


\bibitem[\protect\citeauthoryear{Betcke, Higham, Mehrmann, Schr\"oder, and
  Tisseur}{Betcke et~al\mbox{.}}{2010}]%
        {Betcke2010}
\bibfield{author}{\bibinfo{person}{T. Betcke}, \bibinfo{person}{N.~J. Higham},
  \bibinfo{person}{V. Mehrmann}, \bibinfo{person}{C. Schr\"oder}, {and}
  \bibinfo{person}{F. Tisseur}.} \bibinfo{year}{2010}\natexlab{}.
\newblock \bibinfo{booktitle}{\emph{{NLEVP}: A Collection of Nonlinear
  Eigenvalue Problems}}.
\newblock \bibinfo{type}{{T}echnical {R}eport}.
  \bibinfo{institution}{University of Manchester}.
\newblock


\bibitem[\protect\citeauthoryear{Betcke, Higham, Mehrmann, Schröder, and
  Tisseur}{Betcke et~al\mbox{.}}{2013}]%
        {Betcke:2013:NLEVP}
\bibfield{author}{\bibinfo{person}{Timo Betcke}, \bibinfo{person}{Nicholas~J.
  Higham}, \bibinfo{person}{Volker Mehrmann}, \bibinfo{person}{Christian
  Schröder}, {and} \bibinfo{person}{Fran\c{c}oise Tisseur}.}
  \bibinfo{year}{2013}\natexlab{}.
\newblock \showarticletitle{{NLEVP:} A Collection of Nonlinear Eigenvalue
  Problems}.
\newblock \bibinfo{journal}{\emph{ACM Trans. Math. Softw.}}
  \bibinfo{volume}{39}, \bibinfo{number}{2}, Article \bibinfo{articleno}{7}
  (\bibinfo{date}{Feb.} \bibinfo{year}{2013}), \bibinfo{numpages}{28}~pages.
\newblock
\showISSN{0098-3500}
\urldef\tempurl%
\url{https://doi.org/10.1145/2427023.2427024}
\showDOI{\tempurl}


\bibitem[\protect\citeauthoryear{Betcke and Voss}{Betcke and Voss}{2004}]%
        {Betcke:2004:JD}
\bibfield{author}{\bibinfo{person}{T. Betcke} {and} \bibinfo{person}{H. Voss}.}
  \bibinfo{year}{2004}\natexlab{}.
\newblock \showarticletitle{A {J}acobi-{D}avidson type projection method for
  nonlinear eigenvalue problems}.
\newblock \bibinfo{journal}{\emph{Future Generation Computer Systems}}
  \bibinfo{volume}{20}, \bibinfo{number}{3} (\bibinfo{year}{2004}),
  \bibinfo{pages}{363--372}.
\newblock


\bibitem[\protect\citeauthoryear{Beyn}{Beyn}{2012}]%
        {Beyn:2011:INTEGRAL}
\bibfield{author}{\bibinfo{person}{W.-J. Beyn}.}
  \bibinfo{year}{2012}\natexlab{}.
\newblock \showarticletitle{An integral method for solving nonlinear eigenvalue
  problems}.
\newblock \bibinfo{journal}{\emph{Linear Algebra Appl.}} \bibinfo{volume}{436},
  \bibinfo{number}{10} (\bibinfo{year}{2012}), \bibinfo{pages}{3839--3863}.
\newblock


\bibitem[\protect\citeauthoryear{Beyn, Effenberger, and Kressner}{Beyn
  et~al\mbox{.}}{2011}]%
        {Beyn2011}
\bibfield{author}{\bibinfo{person}{Wolf-J{\"u}rgen Beyn},
  \bibinfo{person}{Cedric Effenberger}, {and} \bibinfo{person}{Daniel
  Kressner}.} \bibinfo{year}{2011}\natexlab{}.
\newblock \showarticletitle{Continuation of eigenvalues and invariant pairs for
  parameterized nonlinear eigenvalue problems}.
\newblock \bibinfo{journal}{\emph{Numer. Math.}} \bibinfo{volume}{119},
  \bibinfo{number}{3} (\bibinfo{date}{10 Jul} \bibinfo{year}{2011}),
  \bibinfo{pages}{489}.
\newblock
\showISSN{0945-3245}
\urldef\tempurl%
\url{https://doi.org/10.1007/s00211-011-0392-1}
\showDOI{\tempurl}


\bibitem[\protect\citeauthoryear{{Beyn} and {Thümmler}}{{Beyn} and
  {Thümmler}}{2010}]%
        {Beyn2010}
\bibfield{author}{\bibinfo{person}{Wolf-Jürgen {Beyn}} {and}
  \bibinfo{person}{Vera {Thümmler}}.} \bibinfo{year}{2010}\natexlab{}.
\newblock \showarticletitle{{Continuation of invariant subspaces for
  parameterized quadratic eigenvalue problems.}}
\newblock \bibinfo{journal}{\emph{SIAM J. Matrix Anal. Appl.}}
  \bibinfo{volume}{31}, \bibinfo{number}{3} (\bibinfo{year}{2010}),
  \bibinfo{pages}{1361--1381}.
\newblock
\showISSN{0895-4798; 1095-7162/e}
\urldef\tempurl%
\url{https://doi.org/10.1137/080723107}
\showDOI{\tempurl}


\bibitem[\protect\citeauthoryear{Bezanson, Edelman, Karpinski, and
  Shah}{Bezanson et~al\mbox{.}}{2017}]%
        {Bezanson2017}
\bibfield{author}{\bibinfo{person}{Jeff Bezanson}, \bibinfo{person}{Alan
  Edelman}, \bibinfo{person}{Stefan Karpinski}, {and} \bibinfo{person}{Viral~B.
  Shah}.} \bibinfo{year}{2017}\natexlab{}.
\newblock \showarticletitle{{Julia}: A Fresh Approach to Numerical Computing}.
\newblock \bibinfo{journal}{\emph{SIAM Rev.}} \bibinfo{volume}{59},
  \bibinfo{number}{1} (\bibinfo{year}{2017}), \bibinfo{pages}{65--98}.
\newblock
\urldef\tempurl%
\url{https://doi.org/10.1137/141000671}
\showDOI{\tempurl}
\showeprint{https://doi.org/10.1137/141000671}


\bibitem[\protect\citeauthoryear{Bueler}{Bueler}{2007}]%
        {Bueler:2007:ERROR}
\bibfield{author}{\bibinfo{person}{E. Bueler}.}
  \bibinfo{year}{2007}\natexlab{}.
\newblock \showarticletitle{Error bounds for approximate eigenvalues of
  periodic-coefficient linear delay differential equations}.
\newblock \bibinfo{journal}{\emph{SIAM J. Numer. Anal.}} \bibinfo{volume}{45},
  \bibinfo{number}{6} (\bibinfo{year}{2007}), \bibinfo{pages}{2510--2536}.
\newblock


\bibitem[\protect\citeauthoryear{E.~Jarlebring}{E.~Jarlebring}{2017}]%
        {Jarlebring:2017:TIAR}
\bibfield{author}{\bibinfo{person}{O.~Runborg E.~Jarlebring, G.~Mele}.}
  \bibinfo{year}{2017}\natexlab{}.
\newblock \showarticletitle{The waveguide eigenvlaue problem and the tensor
  infinite {Arnoldi} method}.
\newblock \bibinfo{journal}{\emph{SIAM J. Sci. Comput.}}  \bibinfo{volume}{39}
  (\bibinfo{year}{2017}), \bibinfo{pages}{A1062--A1088}.
\newblock


\bibitem[\protect\citeauthoryear{Effenberger}{Effenberger}{2013a}]%
        {Effenberger2013a}
\bibfield{author}{\bibinfo{person}{C. Effenberger}.}
  \bibinfo{year}{2013}\natexlab{a}.
\newblock \emph{\bibinfo{title}{Robust Solution Methods for Nonlinear
  Eigenvalue Problems}}.
\newblock \bibinfo{thesistype}{Ph.D. Dissertation}. \bibinfo{school}{EPF
  Lausanne}.
\newblock


\bibitem[\protect\citeauthoryear{Effenberger}{Effenberger}{2013b}]%
        {Effenberger2013}
\bibfield{author}{\bibinfo{person}{C. Effenberger}.}
  \bibinfo{year}{2013}\natexlab{b}.
\newblock \showarticletitle{{Robust successive computation of eigenpairs for
  nonlinear eigenvalue problems}}.
\newblock \bibinfo{journal}{\emph{SIAM J. Matrix Anal. Appl.}}
  \bibinfo{volume}{34}, \bibinfo{number}{3} (\bibinfo{year}{2013}),
  \bibinfo{pages}{1231--1256}.
\newblock
\showISSN{0895-4798; 1095-7162/e}
\urldef\tempurl%
\url{https://doi.org/10.1137/120885644}
\showDOI{\tempurl}


\bibitem[\protect\citeauthoryear{Effenberger and Kressner}{Effenberger and
  Kressner}{2012}]%
        {Effenberger2012a}
\bibfield{author}{\bibinfo{person}{C. Effenberger} {and} \bibinfo{person}{D.
  Kressner}.} \bibinfo{year}{2012}\natexlab{}.
\newblock \showarticletitle{Chebyshev interpolation for nonlinear eigenvalue
  problems}.
\newblock \bibinfo{journal}{\emph{BIT}} \bibinfo{volume}{52},
  \bibinfo{number}{4} (\bibinfo{year}{2012}), \bibinfo{pages}{933--951}.
\newblock


\bibitem[\protect\citeauthoryear{Effenberger, Kressner, and
  Engström}{Effenberger et~al\mbox{.}}{2012}]%
        {Effenberger2012}
\bibfield{author}{\bibinfo{person}{C. Effenberger}, \bibinfo{person}{D.
  Kressner}, {and} \bibinfo{person}{C. Engström}.}
  \bibinfo{year}{2012}\natexlab{}.
\newblock \showarticletitle{Linearization techniques for band structure
  calculations in absorbing photonic crystals}.
\newblock \bibinfo{journal}{\emph{Int. J. Numer. Methods Eng.}}
  \bibinfo{volume}{89}, \bibinfo{number}{2} (\bibinfo{year}{2012}),
  \bibinfo{pages}{180--191}.
\newblock


\bibitem[\protect\citeauthoryear{Engelborghs, Luzyanina, and Roose}{Engelborghs
  et~al\mbox{.}}{2002}]%
        {Engelborghs:2002:DDEBIFTOOL}
\bibfield{author}{\bibinfo{person}{K. Engelborghs}, \bibinfo{person}{T.
  Luzyanina}, {and} \bibinfo{person}{D. Roose}.}
  \bibinfo{year}{2002}\natexlab{}.
\newblock \showarticletitle{Numerical Bifurcation Analysis of Delay
  Differential Equations Using {DDE-BIFTOOL}}.
\newblock \bibinfo{journal}{\emph{ACM Trans. Math. Softw.}}
  \bibinfo{volume}{28}, \bibinfo{number}{1} (\bibinfo{year}{2002}),
  \bibinfo{pages}{1--24}.
\newblock


\bibitem[\protect\citeauthoryear{Engelborghs, Luzyanina, and
  Samaey}{Engelborghs et~al\mbox{.}}{2001}]%
        {Engelborghs:2001:DDEBIFTOOL}
\bibfield{author}{\bibinfo{person}{K. Engelborghs}, \bibinfo{person}{T.
  Luzyanina}, {and} \bibinfo{person}{G. Samaey}.}
  \bibinfo{year}{2001}\natexlab{}.
\newblock \bibinfo{booktitle}{\emph{{DDE-BIFTOOL} v. 2.00: a {M}atlab package
  for bifurcation analysis of delay differential equations}}.
\newblock \bibinfo{type}{{T}echnical {R}eport}.
  \bibinfo{institution}{K.U.Leuven, Leuven, Belgium}.
\newblock


\bibitem[\protect\citeauthoryear{Engström}{Engström}{2010}]%
        {Engstrm2010}
\bibfield{author}{\bibinfo{person}{C. Engström}.}
  \bibinfo{year}{2010}\natexlab{}.
\newblock \showarticletitle{On the spectrum of a holomorphic operator-valued
  function with applications to absorptive photonic crystals}.
\newblock \bibinfo{journal}{\emph{Math. Models Methods Appl. Sci.}}
  \bibinfo{volume}{20}, \bibinfo{number}{8} (\bibinfo{year}{2010}),
  \bibinfo{pages}{1319--1341}.
\newblock


\bibitem[\protect\citeauthoryear{Fliss}{Fliss}{2013}]%
        {Fliss:2017:DTN}
\bibfield{author}{\bibinfo{person}{Sonia Fliss}.}
  \bibinfo{year}{2013}\natexlab{}.
\newblock \showarticletitle{A Dirichlet-to-Neumann Approach for The Exact
  Computation of Guided Modes in Photonic Crystal Waveguides}.
\newblock \bibinfo{journal}{\emph{SIAM J. Sci. Comput.}} \bibinfo{volume}{35},
  \bibinfo{number}{2} (\bibinfo{year}{2013}), \bibinfo{pages}{B438--B461}.
\newblock


\bibitem[\protect\citeauthoryear{Gaaf and Jarlebring}{Gaaf and
  Jarlebring}{2017}]%
        {Gaaf:2017:INFBILANCZOS}
\bibfield{author}{\bibinfo{person}{S.~W. Gaaf} {and} \bibinfo{person}{E.
  Jarlebring}.} \bibinfo{year}{2017}\natexlab{}.
\newblock \showarticletitle{{The infinite bi-Lanczos method for nonlinear
  eigenvalue problems}}.
\newblock \bibinfo{journal}{\emph{SIAM J. Sci. Comput.}} \bibinfo{volume}{39},
  \bibinfo{number}{SIAM J. Sci. Comput.} (\bibinfo{year}{2017}),
  \bibinfo{pages}{S898--S919}.
\newblock
\urldef\tempurl%
\url{https://doi.org/10.1137/16M1084195}
\showDOI{\tempurl}
\showeprint{https://doi.org/10.1137/16M1084195}


\bibitem[\protect\citeauthoryear{Garrett, Bai, and Li}{Garrett
  et~al\mbox{.}}{2016}]%
        {Garrett2016}
\bibfield{author}{\bibinfo{person}{C.~Kristopher Garrett},
  \bibinfo{person}{Zhaojun Bai}, {and} \bibinfo{person}{Ren-Cang Li}.}
  \bibinfo{year}{2016}\natexlab{}.
\newblock \showarticletitle{A nonlinear QR algorithm for banded nonlinear
  eigenvalue problems}.
\newblock \bibinfo{journal}{\emph{ACM Trans. Math. Softw.}}
  \bibinfo{volume}{43}, \bibinfo{number}{1}, Article \bibinfo{articleno}{4}
  (\bibinfo{date}{Aug.} \bibinfo{year}{2016}), \bibinfo{numpages}{19}~pages.
\newblock
\showISSN{0098-3500}
\urldef\tempurl%
\url{https://doi.org/10.1145/2870628}
\showDOI{\tempurl}


\bibitem[\protect\citeauthoryear{G\"uttel, Beeumen, Meerbergen, and
  Michiels}{G\"uttel et~al\mbox{.}}{2014}]%
        {Guttel:2014:NLEIGS}
\bibfield{author}{\bibinfo{person}{S. G\"uttel}, \bibinfo{person}{R.~Van
  Beeumen}, \bibinfo{person}{K. Meerbergen}, {and} \bibinfo{person}{W.
  Michiels}.} \bibinfo{year}{2014}\natexlab{}.
\newblock \showarticletitle{{NLEIGS:} A Class of Fully Rational {Krylov}
  Methods for Nonlinear Eigenvalue Problems}.
\newblock \bibinfo{journal}{\emph{SIAM J. Sci. Comput.}} \bibinfo{volume}{36},
  \bibinfo{number}{6} (\bibinfo{year}{2014}), \bibinfo{pages}{A2842--A2864}.
\newblock


\bibitem[\protect\citeauthoryear{Güttel and Tisseur}{Güttel and
  Tisseur}{2017}]%
        {Guttel2017}
\bibfield{author}{\bibinfo{person}{Stefan Güttel} {and}
  \bibinfo{person}{Francoise Tisseur}.} \bibinfo{year}{2017}\natexlab{}.
\newblock \showarticletitle{The nonlinear eigenvalue problem}.
\newblock \bibinfo{journal}{\emph{Acta Numerica}}  \bibinfo{volume}{26}
  (\bibinfo{year}{2017}), \bibinfo{pages}{1--94}.
\newblock
\urldef\tempurl%
\url{https://doi.org/10.1017/S0962492917000034}
\showDOI{\tempurl}


\bibitem[\protect\citeauthoryear{Hernandez, Roman, and Vidal}{Hernandez
  et~al\mbox{.}}{2003}]%
        {Hernandez:2003:SSL}
\bibfield{author}{\bibinfo{person}{V. Hernandez}, \bibinfo{person}{J.~E.
  Roman}, {and} \bibinfo{person}{V. Vidal}.} \bibinfo{year}{2003}\natexlab{}.
\newblock \showarticletitle{{SLEPc}: {S}calable {L}ibrary for {E}igenvalue
  {P}roblem {C}omputations}.
\newblock \bibinfo{journal}{\emph{Lect. Notes Comput. Sci.}}
  \bibinfo{volume}{2565} (\bibinfo{year}{2003}), \bibinfo{pages}{377--391}.
\newblock


\bibitem[\protect\citeauthoryear{Hernandez, Roman, and Vidal}{Hernandez
  et~al\mbox{.}}{2005}]%
        {Hernandez:2005:SSF}
\bibfield{author}{\bibinfo{person}{Vicente Hernandez}, \bibinfo{person}{Jose~E.
  Roman}, {and} \bibinfo{person}{Vicente Vidal}.}
  \bibinfo{year}{2005}\natexlab{}.
\newblock \showarticletitle{{SLEPc}: A scalable and flexible toolkit for the
  solution of eigenvalue problems}.
\newblock \bibinfo{journal}{\emph{{ACM} Trans. Math. Software}}
  \bibinfo{volume}{31}, \bibinfo{number}{3} (\bibinfo{year}{2005}),
  \bibinfo{pages}{351--362}.
\newblock


\bibitem[\protect\citeauthoryear{Higham}{Higham}{2006}]%
        {Higham:2006:FUNCTIONSOFMATRICES}
\bibfield{author}{\bibinfo{person}{N.~J. Higham}.}
  \bibinfo{year}{2006}\natexlab{}.
\newblock \showarticletitle{Functions of Matrices}.
\newblock In \bibinfo{booktitle}{\emph{Handbook of linear algebra}},
  \bibfield{editor}{\bibinfo{person}{L.~Hogben}} (Ed.). \bibinfo{publisher}{CRC
  Press}.
\newblock


\bibitem[\protect\citeauthoryear{Higham, Li, and Tisseur}{Higham
  et~al\mbox{.}}{2008}]%
        {Higham:2008:BACKWARD}
\bibfield{author}{\bibinfo{person}{N.~J. Higham}, \bibinfo{person}{R.-C. Li},
  {and} \bibinfo{person}{F. Tisseur}.} \bibinfo{year}{2008}\natexlab{}.
\newblock \showarticletitle{Backward error of polynomial eigenproblems solved
  by linearization}.
\newblock \bibinfo{journal}{\emph{SIAM J. Matrix Anal. Appl.}}
  \bibinfo{volume}{29}, \bibinfo{number}{4} (\bibinfo{year}{2008}),
  \bibinfo{pages}{1218–--1241}.
\newblock


\bibitem[\protect\citeauthoryear{Insperger and St\'ep\'an}{Insperger and
  St\'ep\'an}{2002a}]%
        {Insperger:2002:SEMIDISC}
\bibfield{author}{\bibinfo{person}{T. Insperger} {and} \bibinfo{person}{G.
  St\'ep\'an}.} \bibinfo{year}{2002}\natexlab{a}.
\newblock \showarticletitle{Semi-discretization method for delayed systems}.
\newblock \bibinfo{journal}{\emph{Int. J. Numer. Methods Eng.}}
  \bibinfo{volume}{55}, \bibinfo{number}{5} (\bibinfo{year}{2002}),
  \bibinfo{pages}{503--518}.
\newblock


\bibitem[\protect\citeauthoryear{Insperger and St\'ep\'an}{Insperger and
  St\'ep\'an}{2002b}]%
        {Insberger:2002:MATHIEU}
\bibfield{author}{\bibinfo{person}{T. Insperger} {and} \bibinfo{person}{G.
  St\'ep\'an}.} \bibinfo{year}{2002}\natexlab{b}.
\newblock \showarticletitle{Stability chart for the delayed {M}athieu
  equation}.
\newblock \bibinfo{journal}{\emph{Proc. R. Soc. Lond., Ser. A, Math. Phys. Eng.
  Sci.}} \bibinfo{volume}{458}, \bibinfo{number}{2024} (\bibinfo{year}{2002}),
  \bibinfo{pages}{1989--1998}.
\newblock


\bibitem[\protect\citeauthoryear{Jarlebring}{Jarlebring}{2008}]%
        {Jarlebring:2008:THESIS}
\bibfield{author}{\bibinfo{person}{E. Jarlebring}.}
  \bibinfo{year}{2008}\natexlab{}.
\newblock \emph{\bibinfo{title}{The spectrum of delay-differential equations:
  numerical methods, stability and perturbation}}.
\newblock \bibinfo{thesistype}{Ph.D. Dissertation}. \bibinfo{school}{TU
  Braunschweig}.
\newblock


\bibitem[\protect\citeauthoryear{Jarlebring}{Jarlebring}{2017}]%
        {Jarlebring:2017:BROYDEN}
\bibfield{author}{\bibinfo{person}{E. Jarlebring}.}
  \bibinfo{year}{2017}\natexlab{}.
\newblock \bibinfo{booktitle}{\emph{Broyden's method for nonlinear
  eigenproblems}}.
\newblock \bibinfo{type}{{T}echnical {R}eport}. \bibinfo{institution}{KTH Royal
  Institute of Technology}.
\newblock


\bibitem[\protect\citeauthoryear{Jarlebring, Koskela, and Mele}{Jarlebring
  et~al\mbox{.}}{2018}]%
        {Jarlebring:2018:DISGUISED}
\bibfield{author}{\bibinfo{person}{E. Jarlebring}, \bibinfo{person}{A.
  Koskela}, {and} \bibinfo{person}{G. Mele}.} \bibinfo{year}{2018}\natexlab{}.
\newblock \showarticletitle{Disguised and new quasi-Newton methods for
  nonlinear eigenvalue problems}.
\newblock \bibinfo{journal}{\emph{Numer. Algor.}} \bibinfo{volume}{79},
  \bibinfo{number}{1} (\bibinfo{date}{01 Sep} \bibinfo{year}{2018}),
  \bibinfo{pages}{311--335}.
\newblock
\showISSN{1572-9265}
\urldef\tempurl%
\url{https://doi.org/10.1007/s11075-017-0438-2}
\showDOI{\tempurl}


\bibitem[\protect\citeauthoryear{Jarlebring, Meerbergen, and
  Michiels}{Jarlebring et~al\mbox{.}}{2010a}]%
        {Jarlebring:2010:TAYLORARNOLDI}
\bibfield{author}{\bibinfo{person}{Elias Jarlebring}, \bibinfo{person}{Karl
  Meerbergen}, {and} \bibinfo{person}{Wim Michiels}.}
  \bibinfo{year}{2010}\natexlab{a}.
\newblock \showarticletitle{An Arnoldi method with structured starting vectors
  for the delay eigenvalue problem}, In \bibinfo{booktitle}{Proceedings of the
  9th {IFAC} workshop on time-delay systems, Prague}.
\newblock \bibinfo{journal}{\emph{IFAC Proceedings Volumes}}
  \bibinfo{volume}{43}, \bibinfo{number}{2}, \bibinfo{pages}{57 -- 62}.
\newblock
\showISSN{1474-6670}
\urldef\tempurl%
\url{http://www.sciencedirect.com/science/article/pii/S1474667016325034}
\showURL{%
\tempurl}
\newblock
\shownote{9th IFAC Workshop on Time Delay Systems.}


\bibitem[\protect\citeauthoryear{Jarlebring, Meerbergen, and
  Michiels}{Jarlebring et~al\mbox{.}}{2010b}]%
        {Jarlebring:2010:DELAYARNOLDI}
\bibfield{author}{\bibinfo{person}{E. Jarlebring}, \bibinfo{person}{K.
  Meerbergen}, {and} \bibinfo{person}{W. Michiels}.}
  \bibinfo{year}{2010}\natexlab{b}.
\newblock \showarticletitle{A {K}rylov method for the delay eigenvalue
  problem}.
\newblock \bibinfo{journal}{\emph{SIAM J. Sci. Comput.}} \bibinfo{volume}{32},
  \bibinfo{number}{6} (\bibinfo{year}{2010}), \bibinfo{pages}{3278--3300}.
\newblock


\bibitem[\protect\citeauthoryear{Jarlebring, Meerbergen, and
  Michiels}{Jarlebring et~al\mbox{.}}{2014}]%
        {Jarlebring:2017:SCHUR}
\bibfield{author}{\bibinfo{person}{E. Jarlebring}, \bibinfo{person}{K.
  Meerbergen}, {and} \bibinfo{person}{W. Michiels}.}
  \bibinfo{year}{2014}\natexlab{}.
\newblock \showarticletitle{Computing a partial {Schur} factorization of
  nonlinear eigenvalue problems using the infinite {Arnoldi} method}.
\newblock \bibinfo{journal}{\emph{SIAM J. Matrix Anal. Appl.}}
  \bibinfo{volume}{35}, \bibinfo{number}{2} (\bibinfo{year}{2014}),
  \bibinfo{pages}{411--436}.
\newblock
\urldef\tempurl%
\url{https://doi.org/10.1137/110858148}
\showDOI{\tempurl}
\showeprint{https://doi.org/10.1137/110858148}


\bibitem[\protect\citeauthoryear{Jarlebring, Michiels, and
  Meerbergen}{Jarlebring et~al\mbox{.}}{2012}]%
        {Jarlebring:2012:INFARNOLDI}
\bibfield{author}{\bibinfo{person}{E. Jarlebring}, \bibinfo{person}{W.
  Michiels}, {and} \bibinfo{person}{K. Meerbergen}.}
  \bibinfo{year}{2012}\natexlab{}.
\newblock \showarticletitle{A linear eigenvalue algorithm for the nonlinear
  eigenvalue problem}.
\newblock \bibinfo{journal}{\emph{Numer. Math.}} \bibinfo{volume}{122},
  \bibinfo{number}{1} (\bibinfo{year}{2012}), \bibinfo{pages}{169--195}.
\newblock


\bibitem[\protect\citeauthoryear{Jarlebring and Voss}{Jarlebring and
  Voss}{2005}]%
        {Jarlebring:2005:RATIONAL}
\bibfield{author}{\bibinfo{person}{E. Jarlebring} {and} \bibinfo{person}{H.
  Voss}.} \bibinfo{year}{2005}\natexlab{}.
\newblock \showarticletitle{Rational Krylov for nonlinear eigenproblems, an
  iterative projection method}.
\newblock \bibinfo{journal}{\emph{Appl. Math.}} \bibinfo{volume}{50},
  \bibinfo{number}{6} (\bibinfo{year}{2005}), \bibinfo{pages}{543--554}.
\newblock


\bibitem[\protect\citeauthoryear{Kaufman}{Kaufman}{2006}]%
        {Kaufman2006}
\bibfield{author}{\bibinfo{person}{L. Kaufman}.}
  \bibinfo{year}{2006}\natexlab{}.
\newblock \showarticletitle{Eigenvalue problems in fiber optic design}.
\newblock \bibinfo{journal}{\emph{SIAM J. Matrix Anal. Appl.}}
  \bibinfo{volume}{28}, \bibinfo{number}{1} (\bibinfo{year}{2006}),
  \bibinfo{pages}{105--117}.
\newblock


\bibitem[\protect\citeauthoryear{Kressner}{Kressner}{2009}]%
        {Kressner:2009:BLOCKNEWTON}
\bibfield{author}{\bibinfo{person}{D. Kressner}.}
  \bibinfo{year}{2009}\natexlab{}.
\newblock \showarticletitle{A block {N}ewton method for nonlinear eigenvalue
  problems}.
\newblock \bibinfo{journal}{\emph{Numer. Math.}} \bibinfo{volume}{114},
  \bibinfo{number}{2} (\bibinfo{year}{2009}), \bibinfo{pages}{355--372}.
\newblock


\bibitem[\protect\citeauthoryear{Kublanovskaya}{Kublanovskaya}{1970}]%
        {Kublanovskaya:1970:APPROACH}
\bibfield{author}{\bibinfo{person}{V. Kublanovskaya}.}
  \bibinfo{year}{1970}\natexlab{}.
\newblock \showarticletitle{On an approach to the solution of the generalized
  latent value problem for $\lambda$-matrices}.
\newblock \bibinfo{journal}{\emph{SIAM J. Numer. Anal.}}  \bibinfo{volume}{7}
  (\bibinfo{year}{1970}), \bibinfo{pages}{532--537}.
\newblock


\bibitem[\protect\citeauthoryear{Mehrmann and Voss}{Mehrmann and Voss}{2004}]%
        {Mehrmann:2004:NLEVP}
\bibfield{author}{\bibinfo{person}{V. Mehrmann} {and} \bibinfo{person}{H.
  Voss}.} \bibinfo{year}{2004}\natexlab{}.
\newblock \showarticletitle{Nonlinear eigenvalue problems: A Challange for
  modern eigenvalue methods}.
\newblock \bibinfo{journal}{\emph{GAMM Mitteilungen}}  \bibinfo{volume}{27}
  (\bibinfo{year}{2004}), \bibinfo{pages}{121--152}.
\newblock


\bibitem[\protect\citeauthoryear{Mele and Jarlebring}{Mele and
  Jarlebring}{2018}]%
        {Mele2018}
\bibfield{author}{\bibinfo{person}{Giampaolo Mele} {and} \bibinfo{person}{Elias
  Jarlebring}.} \bibinfo{year}{2018}\natexlab{}.
\newblock \showarticletitle{On restarting the tensor infinite {Arnoldi}
  method}.
\newblock \bibinfo{journal}{\emph{BIT}} \bibinfo{volume}{58},
  \bibinfo{number}{1} (\bibinfo{date}{01 Mar} \bibinfo{year}{2018}),
  \bibinfo{pages}{133--162}.
\newblock
\showISSN{1572-9125}
\urldef\tempurl%
\url{https://doi.org/10.1007/s10543-017-0671-z}
\showDOI{\tempurl}


\bibitem[\protect\citeauthoryear{Michiels and Niculescu}{Michiels and
  Niculescu}{2007}]%
        {Michiels:2007:STABILITYBOOK}
\bibfield{author}{\bibinfo{person}{W. Michiels} {and} \bibinfo{person}{S.-I.
  Niculescu}.} \bibinfo{year}{2007}\natexlab{}.
\newblock \bibinfo{booktitle}{\emph{Stability and Stabilization of Time-Delay
  Systems: An Eigenvalue-Based Approach}}.
\newblock \bibinfo{publisher}{SIAM Publications, Philadelphia}.
\newblock


\bibitem[\protect\citeauthoryear{Neumaier}{Neumaier}{1985}]%
        {Neumaier:1985:RESINV}
\bibfield{author}{\bibinfo{person}{A. Neumaier}.}
  \bibinfo{year}{1985}\natexlab{}.
\newblock \showarticletitle{Residual inverse iteration for the nonlinear
  eigenvalue problem}.
\newblock \bibinfo{journal}{\emph{SIAM J. Numer. Anal.}}  \bibinfo{volume}{22}
  (\bibinfo{year}{1985}), \bibinfo{pages}{914--923}.
\newblock


\bibitem[\protect\citeauthoryear{Ooi, Mizuno, Sogabe, Yamamoto, and Zhang}{Ooi
  et~al\mbox{.}}{2017}]%
        {Ooi:2017:NEP}
\bibfield{author}{\bibinfo{person}{Kouhei Ooi}, \bibinfo{person}{Yoshinori
  Mizuno}, \bibinfo{person}{Tomohiro Sogabe}, \bibinfo{person}{Yusaku
  Yamamoto}, {and} \bibinfo{person}{Shao-Liang Zhang}.}
  \bibinfo{year}{2017}\natexlab{}.
\newblock \showarticletitle{Solution of a nonlinear eigenvalue problem using
  signed singular values}.
\newblock \bibinfo{journal}{\emph{East Asian Journal on Applied Mathematics}}
  \bibinfo{volume}{7}, \bibinfo{number}{4} (\bibinfo{year}{2017}),
  \bibinfo{pages}{799–809}.
\newblock
\urldef\tempurl%
\url{https://doi.org/10.4208/eajam.181016.300517c}
\showDOI{\tempurl}


\bibitem[\protect\citeauthoryear{Polizzi}{Polizzi}{2009}]%
        {Polizzi:2009:FEAST}
\bibfield{author}{\bibinfo{person}{Eric Polizzi}.}
  \bibinfo{year}{2009}\natexlab{}.
\newblock \showarticletitle{Density-matrix-based algorithm for solving
  eigenvalue problems}.
\newblock \bibinfo{journal}{\emph{Phys. Rev. B}}  \bibinfo{volume}{79}
  (\bibinfo{date}{Mar} \bibinfo{year}{2009}), \bibinfo{pages}{115112}.
\newblock
Issue 11.
\urldef\tempurl%
\url{https://doi.org/10.1103/PhysRevB.79.115112}
\showDOI{\tempurl}


\bibitem[\protect\citeauthoryear{Regier, Pamnany, Fischer, Noack, Lam, Revels,
  Howard, Giordano, Schlegel, McAuliffe, and Thomas}{Regier
  et~al\mbox{.}}{2018}]%
        {celeste2018}
\bibfield{author}{\bibinfo{person}{Jeffrey Regier}, \bibinfo{person}{Kiran
  Pamnany}, \bibinfo{person}{Keno Fischer}, \bibinfo{person}{Andreas Noack},
  \bibinfo{person}{Maximilian Lam}, \bibinfo{person}{Jarrett Revels},
  \bibinfo{person}{Steve Howard}, \bibinfo{person}{Ryan Giordano},
  \bibinfo{person}{David Schlegel}, \bibinfo{person}{Jon McAuliffe}, {and}
  \bibinfo{person}{Rollin Thomas}.} \bibinfo{year}{2018}\natexlab{}.
\newblock \bibinfo{booktitle}{\emph{Cataloging the Visible Universe through
  Bayesian Inference at Petascale}}.
\newblock \bibinfo{type}{{T}echnical {R}eport}.
  \bibinfo{institution}{University of California, Berkeley}.
\newblock
\newblock
\shownote{arXiv preprint arXiv:1801.10277.}


\bibitem[\protect\citeauthoryear{Ringh, Mele, Karlsson, and Jarlebring}{Ringh
  et~al\mbox{.}}{2018}]%
        {Ringh:2018:SYLVPRECOND}
\bibfield{author}{\bibinfo{person}{E. Ringh}, \bibinfo{person}{G. Mele},
  \bibinfo{person}{J. Karlsson}, {and} \bibinfo{person}{E. Jarlebring}.}
  \bibinfo{year}{2018}\natexlab{}.
\newblock \showarticletitle{Sylvester-based preconditioning for the waveguide
  eigenvalue problem}.
\newblock \bibinfo{journal}{\emph{Linear Algebra Appl.}}  \bibinfo{volume}{542}
  (\bibinfo{year}{2018}), \bibinfo{pages}{441 -- 463}.
\newblock
\showISSN{0024-3795}
\urldef\tempurl%
\url{https://doi.org/10.1016/j.laa.2017.06.027}
\showDOI{\tempurl}
\newblock
\shownote{Proceedings of the 20th ILAS Conference, Leuven, Belgium 2016.}


\bibitem[\protect\citeauthoryear{Roman, Campos, Romero, and Tomas}{Roman
  et~al\mbox{.}}{2018}]%
        {Roman:2018:SLEPC}
\bibfield{author}{\bibinfo{person}{J.~E. Roman}, \bibinfo{person}{C. Campos},
  \bibinfo{person}{E. Romero}, {and} \bibinfo{person}{A. Tomas}.}
  \bibinfo{year}{2018}\natexlab{}.
\newblock \bibinfo{booktitle}{\emph{{SLEPc} Users Manual}}.
\newblock \bibinfo{type}{{T}echnical {R}eport} DSIC-II/24/02 - Revision 3.10.
  \bibinfo{institution}{D. Sistemes Inform\`atics i Computaci\'o, Universitat
  Polit\`ecnica de Val\`encia}.
\newblock


\bibitem[\protect\citeauthoryear{Rott and Jarlebring}{Rott and
  Jarlebring}{2010}]%
        {Rott:2010:ITERATIVE}
\bibfield{author}{\bibinfo{person}{O. Rott} {and} \bibinfo{person}{E.
  Jarlebring}.} \bibinfo{year}{2010}\natexlab{}.
\newblock \showarticletitle{An iterative method for the multipliers of periodic
  delay-differential equations and the analysis of a {PDE} milling model}. In
  \bibinfo{booktitle}{\emph{Proceedings of the 9th {IFAC} workshop on
  time-delay systems, Prague}}. \bibinfo{pages}{1--6}.
\newblock


\bibitem[\protect\citeauthoryear{Ruhe}{Ruhe}{1973}]%
        {Ruhe:1973:NLEVP}
\bibfield{author}{\bibinfo{person}{A. Ruhe}.} \bibinfo{year}{1973}\natexlab{}.
\newblock \showarticletitle{Algorithms for the nonlinear eigenvalue problem}.
\newblock \bibinfo{journal}{\emph{SIAM J. Numer. Anal.}}  \bibinfo{volume}{10}
  (\bibinfo{year}{1973}), \bibinfo{pages}{674--689}.
\newblock


\bibitem[\protect\citeauthoryear{Schmidt and Kauf}{Schmidt and Kauf}{2009}]%
        {SCHMIDT:2009:BAND}
\bibfield{author}{\bibinfo{person}{K. Schmidt} {and} \bibinfo{person}{P.
  Kauf}.} \bibinfo{year}{2009}\natexlab{}.
\newblock \showarticletitle{Computation of the band structure of
  two-dimensional photonic crystals with hp finite elements}.
\newblock \bibinfo{journal}{\emph{Computer Methods in Applied Mechanics and
  Engineering}} \bibinfo{volume}{198}, \bibinfo{number}{13}
  (\bibinfo{year}{2009}), \bibinfo{pages}{1249 -- 1259}.
\newblock
\showISSN{0045-7825}
\urldef\tempurl%
\url{https://doi.org/10.1016/j.cma.2008.06.009}
\showDOI{\tempurl}
\newblock
\shownote{HOFEM07.}


\bibitem[\protect\citeauthoryear{Schreiber}{Schreiber}{2008}]%
        {Schreiber:2008:PHD}
\bibfield{author}{\bibinfo{person}{K. Schreiber}.}
  \bibinfo{year}{2008}\natexlab{}.
\newblock \emph{\bibinfo{title}{Nonlinear Eigenvalue Problems: {N}ewton-type
  Methods and Nonlinear {R}ayleigh Functionals}}.
\newblock \bibinfo{thesistype}{Ph.D. Dissertation}. \bibinfo{school}{{TU}
  {B}erlin}.
\newblock


\bibitem[\protect\citeauthoryear{Schwetlick and Schreiber}{Schwetlick and
  Schreiber}{2012}]%
        {Schwetlick:2012:NONLINEAR}
\bibfield{author}{\bibinfo{person}{Hubert Schwetlick} {and}
  \bibinfo{person}{Kathrin Schreiber}.} \bibinfo{year}{2012}\natexlab{}.
\newblock \showarticletitle{Nonlinear {Rayleigh} functionals}.
\newblock \bibinfo{journal}{\emph{Linear Algebra Appl.}} \bibinfo{volume}{436},
  \bibinfo{number}{10} (\bibinfo{year}{2012}), \bibinfo{pages}{3991 -- 4016}.
\newblock
\showISSN{0024-3795}
\urldef\tempurl%
\url{http://www.sciencedirect.com/science/article/pii/S0024379510003447}
\showURL{%
\tempurl}
\newblock
\shownote{Special Issue dedicated to Heinrich Voss's 65th birthday.}


\bibitem[\protect\citeauthoryear{Shayer and Campbell}{Shayer and
  Campbell}{2000}]%
        {Shayer:2000:STABILITY}
\bibfield{author}{\bibinfo{person}{L.~P. Shayer} {and} \bibinfo{person}{S.~A.
  Campbell}.} \bibinfo{year}{2000}\natexlab{}.
\newblock \showarticletitle{Stability, bifurcation, and multistability in a
  system of two coupled neurons with multiple time delays}.
\newblock \bibinfo{journal}{\emph{SIAM J. Appl. Math.}} \bibinfo{volume}{61},
  \bibinfo{number}{2} (\bibinfo{year}{2000}), \bibinfo{pages}{673--700}.
\newblock


\bibitem[\protect\citeauthoryear{Sleijpen, Booten, Fokkema, and {v}an~{d}er
  Vorst}{Sleijpen et~al\mbox{.}}{1996}]%
        {Sleijpen1996}
\bibfield{author}{\bibinfo{person}{G. Sleijpen}, \bibinfo{person}{A.~G.
  Booten}, \bibinfo{person}{D.~R. Fokkema}, {and} \bibinfo{person}{H.~A.
  {v}an~{d}er Vorst}.} \bibinfo{year}{1996}\natexlab{}.
\newblock \showarticletitle{Jacobi-{D}avidson type methods for generalized
  eigenproblems and polynomial eigenproblems}.
\newblock \bibinfo{journal}{\emph{BIT}} \bibinfo{volume}{36},
  \bibinfo{number}{3} (\bibinfo{year}{1996}), \bibinfo{pages}{595--633}.
\newblock


\bibitem[\protect\citeauthoryear{\'{S}migaj, Betcke, Arridge, Phillips, and
  Schweiger}{\'{S}migaj et~al\mbox{.}}{2015}]%
        {Smigaj:2015:BEM}
\bibfield{author}{\bibinfo{person}{Wojciech \'{S}migaj}, \bibinfo{person}{Timo
  Betcke}, \bibinfo{person}{Simon Arridge}, \bibinfo{person}{Joel Phillips},
  {and} \bibinfo{person}{Martin Schweiger}.} \bibinfo{year}{2015}\natexlab{}.
\newblock \showarticletitle{Solving Boundary Integral Problems with BEM++}.
\newblock \bibinfo{journal}{\emph{ACM Trans. Math. Softw.}}
  \bibinfo{volume}{41}, Article \bibinfo{articleno}{6} (\bibinfo{date}{Feb.}
  \bibinfo{year}{2015}), \bibinfo{numpages}{40}~pages.
\newblock
\showISSN{0098-3500}
\urldef\tempurl%
\url{https://doi.org/10.1145/2590830}
\showDOI{\tempurl}


\bibitem[\protect\citeauthoryear{Spence and Poulton}{Spence and
  Poulton}{2005}]%
        {Spence:2005:PHOTONIC}
\bibfield{author}{\bibinfo{person}{A. Spence} {and} \bibinfo{person}{C.
  Poulton}.} \bibinfo{year}{2005}\natexlab{}.
\newblock \showarticletitle{Photonic band structure calculations using
  nonlinear eigenvalue techniques}.
\newblock \bibinfo{journal}{\emph{J. Comput. Phys.}} \bibinfo{volume}{204},
  \bibinfo{number}{1} (\bibinfo{year}{2005}), \bibinfo{pages}{65--81}.
\newblock


\bibitem[\protect\citeauthoryear{Steinbach}{Steinbach}{2007}]%
        {Steinbach:2007:BEM}
\bibfield{author}{\bibinfo{person}{O. Steinbach}.}
  \bibinfo{year}{2007}\natexlab{}.
\newblock \bibinfo{booktitle}{\emph{Numerical approximation methods for
  elliptic boundary value problems: finite and boundary elements}}.
\newblock \bibinfo{publisher}{Springer}.
\newblock


\bibitem[\protect\citeauthoryear{Steinbach and Unger}{Steinbach and
  Unger}{2009}]%
        {Steinbach:2009:BEM}
\bibfield{author}{\bibinfo{person}{O. Steinbach} {and} \bibinfo{person}{G.
  Unger}.} \bibinfo{year}{2009}\natexlab{}.
\newblock \showarticletitle{A boundary element method for the {Dirichlet}
  eigenvalue problem of the {Laplace} operator}.
\newblock \bibinfo{journal}{\emph{Numer. Math.}}  \bibinfo{volume}{113}
  (\bibinfo{year}{2009}), \bibinfo{pages}{281--298}.
\newblock


\bibitem[\protect\citeauthoryear{{Szyld} and {Xue}}{{Szyld} and {Xue}}{2013}]%
        {Szyld2013}
\bibfield{author}{\bibinfo{person}{D.~B. {Szyld}} {and} \bibinfo{person}{F.
  {Xue}}.} \bibinfo{year}{2013}\natexlab{}.
\newblock \showarticletitle{Local convergence analysis of several inexact
  {Newton-type algorithms} for general nonlinear eigenvalue problems.}
\newblock \bibinfo{journal}{\emph{Numer. Math.}} \bibinfo{volume}{123},
  \bibinfo{number}{2} (\bibinfo{year}{2013}), \bibinfo{pages}{333--362}.
\newblock
\showISSN{0029-599X; 0945-3245/e}
\urldef\tempurl%
\url{https://doi.org/10.1007/s00211-012-0489-1}
\showDOI{\tempurl}


\bibitem[\protect\citeauthoryear{Szyld and Xue}{Szyld and Xue}{2014}]%
        {Szyld2014}
\bibfield{author}{\bibinfo{person}{D.~B. Szyld} {and} \bibinfo{person}{F.
  Xue}.} \bibinfo{year}{2014}\natexlab{}.
\newblock \showarticletitle{Several properties of invariant pairs of nonlinear
  algebraic eigenvalue problems}.
\newblock \bibinfo{journal}{\emph{IMA J. Numer. Anal.}}  \bibinfo{volume}{34}
  (\bibinfo{year}{2014}), \bibinfo{pages}{921–954}.
\newblock


\bibitem[\protect\citeauthoryear{{Szyld} and {Xue}}{{Szyld} and {Xue}}{2015}]%
        {Szyld:2015:LOCAL}
\bibfield{author}{\bibinfo{person}{D.~B. {Szyld}} {and} \bibinfo{person}{F.
  {Xue}}.} \bibinfo{year}{2015}\natexlab{}.
\newblock \showarticletitle{Local convergence of {Newton-like} methods for
  degenerate eigenvalues of nonlinear eigenproblems}.
\newblock \bibinfo{journal}{\emph{Numer. Math.}} \bibinfo{volume}{129},
  \bibinfo{number}{2} (\bibinfo{year}{2015}), \bibinfo{pages}{353--381}.
\newblock
\showISSN{0029-599X; 0945-3245/e}
\urldef\tempurl%
\url{https://doi.org/10.1007/s00211-014-0639-8}
\showDOI{\tempurl}


\bibitem[\protect\citeauthoryear{Tausch and Butler}{Tausch and Butler}{2000}]%
        {Tausch2000}
\bibfield{author}{\bibinfo{person}{J. Tausch} {and} \bibinfo{person}{J.
  Butler}.} \bibinfo{year}{2000}\natexlab{}.
\newblock \showarticletitle{Floquet multipliers of periodic waveguides via
  {Dirichlet-to-Neumann} maps}.
\newblock \bibinfo{journal}{\emph{J. Comput. Phys.}} \bibinfo{volume}{159},
  \bibinfo{number}{1} (\bibinfo{year}{2000}), \bibinfo{pages}{90--102}.
\newblock


\bibitem[\protect\citeauthoryear{Tausch and Butler}{Tausch and Butler}{2002}]%
        {Tausch2002}
\bibfield{author}{\bibinfo{person}{J. Tausch} {and} \bibinfo{person}{J.
  Butler}.} \bibinfo{year}{2002}\natexlab{}.
\newblock \showarticletitle{Efficient analysis of periodic dielectric
  waveguides using {Dirichlet-to-Neumann} maps}.
\newblock \bibinfo{journal}{\emph{J Opt Soc Am A Opt Image Sci Vis.}}
  \bibinfo{volume}{19}, \bibinfo{number}{6} (\bibinfo{year}{2002}),
  \bibinfo{pages}{1120--8}.
\newblock


\bibitem[\protect\citeauthoryear{Tisseur and Meerbergen}{Tisseur and
  Meerbergen}{2001}]%
        {Tisseur:2001:QUADRATIC}
\bibfield{author}{\bibinfo{person}{F. Tisseur} {and} \bibinfo{person}{K.
  Meerbergen}.} \bibinfo{year}{2001}\natexlab{}.
\newblock \showarticletitle{The quadratic eigenvalue problem}.
\newblock \bibinfo{journal}{\emph{SIAM Rev.}} \bibinfo{volume}{43},
  \bibinfo{number}{2} (\bibinfo{year}{2001}), \bibinfo{pages}{235--286}.
\newblock


\bibitem[\protect\citeauthoryear{Unger}{Unger}{2013}]%
        {unger2013convergence}
\bibfield{author}{\bibinfo{person}{G. Unger}.} \bibinfo{year}{2013}\natexlab{}.
\newblock \showarticletitle{Convergence orders of iterative methods for
  nonlinear eigenvalue problems}.
\newblock In \bibinfo{booktitle}{\emph{Advanced Finite Element Methods and
  Applications}}. \bibinfo{publisher}{Springer}, \bibinfo{pages}{217--237}.
\newblock


\bibitem[\protect\citeauthoryear{Unger}{Unger}{1950}]%
        {Unger:1950:NICHTLINEARE}
\bibfield{author}{\bibinfo{person}{H. Unger}.} \bibinfo{year}{1950}\natexlab{}.
\newblock \showarticletitle{{Nichtlineare Behandlung von Eigenwertaufgaben}}.
\newblock \bibinfo{journal}{\emph{Z. Angew. Math. Mech.}}  \bibinfo{volume}{30}
  (\bibinfo{year}{1950}), \bibinfo{pages}{281--282}.
\newblock
\newblock
\shownote{English translation:
  http://www.math.tu-dresden.de/\~{}schwetli/Unger.html.}


\bibitem[\protect\citeauthoryear{{Van Beeumen}, Meerbergen, and Michiels}{{Van
  Beeumen} et~al\mbox{.}}{2013}]%
        {VanBeeumen:2013:RATIONAL}
\bibfield{author}{\bibinfo{person}{R. {Van Beeumen}}, \bibinfo{person}{K.
  Meerbergen}, {and} \bibinfo{person}{W. Michiels}.}
  \bibinfo{year}{2013}\natexlab{}.
\newblock \showarticletitle{A rational {Krylov} method based on {Hermite}
  interpolation for Nonlinear Eigenvalue Problems}.
\newblock \bibinfo{journal}{\emph{SIAM J. Sci. Comput.}} \bibinfo{volume}{35},
  \bibinfo{number}{1} (\bibinfo{year}{2013}), \bibinfo{pages}{A327--A350}.
\newblock


\bibitem[\protect\citeauthoryear{{Van Beeumen}, Meerbergen, and Michiels}{{Van
  Beeumen} et~al\mbox{.}}{2015}]%
        {VanBeeumen:2015:COMPACTRATIONAL}
\bibfield{author}{\bibinfo{person}{R. {Van Beeumen}}, \bibinfo{person}{K.
  Meerbergen}, {and} \bibinfo{person}{W. Michiels}.}
  \bibinfo{year}{2015}\natexlab{}.
\newblock \showarticletitle{Compact rational {Krylov} methods for Nonlinear
  Eigenvalue Problems}.
\newblock \bibinfo{journal}{\emph{SIAM J. Sci. Comput.}} \bibinfo{volume}{36},
  \bibinfo{number}{2} (\bibinfo{year}{2015}), \bibinfo{pages}{820--838}.
\newblock


\bibitem[\protect\citeauthoryear{Voss}{Voss}{2003}]%
        {Voss:2003:MAXMIN}
\bibfield{author}{\bibinfo{person}{H. Voss}.} \bibinfo{year}{2003}\natexlab{}.
\newblock \showarticletitle{A maxmin principle for nonlinear eigenvalue
  problems with application to a rational spectral problem in fluid-solid
  vibration}.
\newblock \bibinfo{journal}{\emph{Appl. Math., Praha}} \bibinfo{volume}{48},
  \bibinfo{number}{6} (\bibinfo{year}{2003}), \bibinfo{pages}{607--622}.
\newblock


\bibitem[\protect\citeauthoryear{Voss}{Voss}{2004}]%
        {Voss:2004:ARNOLDI}
\bibfield{author}{\bibinfo{person}{H. Voss}.} \bibinfo{year}{2004}\natexlab{}.
\newblock \showarticletitle{An {Arnoldi} method for nonlinear eigenvalue
  problems}.
\newblock \bibinfo{journal}{\emph{BIT}}  \bibinfo{volume}{44}
  (\bibinfo{year}{2004}), \bibinfo{pages}{387 -- 401}.
\newblock


\bibitem[\protect\citeauthoryear{Voss}{Voss}{2005}]%
        {Voss2005}
\bibfield{author}{\bibinfo{person}{H. Voss}.} \bibinfo{year}{2005}\natexlab{}.
\newblock \showarticletitle{Locating real eigenvalues of a spectral problem in
  fluid-solid type structures}.
\newblock \bibinfo{journal}{\emph{J. Appl. Math.}} \bibinfo{volume}{2005},
  \bibinfo{number}{1} (\bibinfo{year}{2005}), \bibinfo{pages}{37--48}.
\newblock


\bibitem[\protect\citeauthoryear{Voss}{Voss}{2007}]%
        {Voss:2007:JD}
\bibfield{author}{\bibinfo{person}{H. Voss}.} \bibinfo{year}{2007}\natexlab{}.
\newblock \showarticletitle{A new justification of the Jacobi-Davidson method
  for large eigenproblems}.
\newblock \bibinfo{journal}{\emph{Linear Algebra Appl.}}  \bibinfo{volume}{424}
  (\bibinfo{year}{2007}), \bibinfo{pages}{448--455}.
\newblock


\bibitem[\protect\citeauthoryear{Voss}{Voss}{2012}]%
        {Voss2012}
\bibfield{author}{\bibinfo{person}{H. Voss}.} \bibinfo{year}{2012}\natexlab{}.
\newblock \bibinfo{booktitle}{\emph{Chapter Nonlinear Eigenvalue Problems}}.
\newblock \bibinfo{publisher}{CRC press}.
\newblock
\newblock
\shownote{Handbook in Linear Algebra.}


\bibitem[\protect\citeauthoryear{Xiao, Zhang, Huang, and Sakurai}{Xiao
  et~al\mbox{.}}{2017}]%
        {Xiao:2017:SOLVING}
\bibfield{author}{\bibinfo{person}{J. Xiao}, \bibinfo{person}{C. Zhang},
  \bibinfo{person}{T.-M. Huang}, {and} \bibinfo{person}{T. Sakurai}.}
  \bibinfo{year}{2017}\natexlab{}.
\newblock \showarticletitle{Solving large-scale nonlinear eigenvalue problems
  by rational interpolation and resolvent sampling based {Rayleigh-Ritz}
  method}.
\newblock \bibinfo{journal}{\emph{Internat. J. Numer. Methods Engrg}}
  \bibinfo{volume}{110}, \bibinfo{number}{8} (\bibinfo{year}{2017}),
  \bibinfo{pages}{776--800}.
\newblock


\bibitem[\protect\citeauthoryear{Xue}{Xue}{2018}]%
        {Xue:2018:BHP}
\bibfield{author}{\bibinfo{person}{F. Xue}.} \bibinfo{year}{2018}\natexlab{}.
\newblock \showarticletitle{A Block Preconditioned Harmonic Projection Method
  for Large-Scale Nonlinear Eigenvalue Problems}.
\newblock \bibinfo{journal}{\emph{SIAM J. Matrix Anal. Appl.}}
  \bibinfo{volume}{40}, \bibinfo{number}{3} (\bibinfo{year}{2018}),
  \bibinfo{pages}{A1809--A1835}.
\newblock
\urldef\tempurl%
\url{https://doi.org/10.1137/17M112141X}
\showDOI{\tempurl}


\end{thebibliography}
                             % Sample .bib file with references that match those in
                             % the 'Specifications Document (V1.5)' as well containing
                             % 'legacy' bibs and bibs with 'alternate codings'.
                             % Gerry Murray - March 2012

% History dates
%\received{February 2007}{March 2009}{June 2009}

% Electronic Appendix
%\elecappendix
\appendix
%\medskip

%\section{This is an example of Appendix section head}

\section{The fiber benchmark}\label{sec:fiber}

The benchmark problem in \cite{Betcke:2013:NLEVP} called ``fiber'' contains a
term defined as
\[
   f(\lambda)=g(\sqrt{\lambda}L)
   \]
 where
\[
  g(x)=\frac{L+0.5}{L^2}x\frac{K'_1(x)}{K_1(x)}
\]
Bessel matrix functions are not available, so we use an interpolation approach
to create a matrix function.
\[
  g(x)=\frac{\alpha(x)}{\beta(x)}
  \]
  where
  %\begin{subequations}
  \begin{eqnarray}
  \alpha(x)&=&\frac{L+0.5}{L^2}\frac{x}{K_1(x)^2}\\
  \beta(x)&=&\frac{1}{K'_1(x)K_1(x)}
  \end{eqnarray}
  %\end{subequations}
  The functions $\alpha$ and $\beta$ are selected such that we can carry out
  polynomial interpolation. Note that $K_1(x)$ has a singularity at zero.
  We create Newton polynomials which interpolate $\alpha$ and $\beta$
  in certain interpolation points. The interpolation is
  carried out in \verb#BigFloat#, and subsequently rounded to \verb#Float64#
  in order to lessen the impact of round-off error with many interpolation
  points.
% Acknowledgments

\section{Interface to iar Chebyshev variant}\label{sec:iar_chebyshev_comp}
The method \verb#iar_chebyshev# requires, at each iteration, the computation of the vector $y_0$ defined in \cite[(22)]{Jarlebring:2012:INFARNOLDI}. We refer to that paper for all the notation we use. Our interface handles this computation for \verb#PEP#, \verb#DEP#, with the derivation of \cite{Jarlebring:2012:INFARNOLDI}, and we have further derived the analogous formula for \verb#SPMF_NEP#. More precisely, by using the Taylor series expansion on \eqref{eq:SPMF_NEP}, it holds
\begin{align*}
 y_0= \sum_{i=0}^{m} A_i
 X
 b_i \left( D_N \right)
 \TTh_N (0)
 - Y \TTh_N (0)
\end{align*}
where $\TTh_N(\theta):=(\hat T_0(\theta), \hat T_1(\theta), \dots, \hat T_{N}(\theta))^T$, $D_N$ is the derivation matrix in Chebyshev basis, defined as
$$D_N:=
\begin{pmatrix}
 0 \\ I_{N,N+1} L_{N+1}^{-1}
\end{pmatrix} \in \mathbb{R}^{(n+1) \times (n+1)}
$$
and derived from \cite[(21)]{Jarlebring:2012:INFARNOLDI} and $b_i(\lambda)=(f_i(0)-f_i(\lambda))/\lambda=f_i[\lambda,0]$. With the same technique, this formula can be extended to the computation of $\tilde y_0$ for the shifted and scaled problem \eqref{eq:M_shift_scale} without explicitly constructing \eqref{eq:M_shift_scale} but directly using \eqref{eq:SPMF_NEP} as follows
\begin{align*}
 \tilde y_0 &=
 -\alpha \sum_{i=0}^{m} M_i
 X
 b_i (\sigma I+\alpha D_N)
 \TTh_N (0)
 - Y \TTh_N (0),
\end{align*}
where $b_i(\sigma + \alpha \lambda)=(f_i(\sigma)-f_i(\sigma+\alpha \lambda))/\lambda=
-\alpha f[\sigma+\alpha \lambda,\sigma]$.
The computation of the divided differences matrices
$b_i(\sigma I + \alpha D_N)=-\alpha f_i \left[\sigma I + \alpha D_N, \sigma I \right]$
can be carried out in terms of functions defining the original problem \eqref{eq:SPMF_NEP} by applying the following result, which is a consequence of the theory for Fréchet derivative in \cite[Section~3.2]{Higham:2006:FUNCTIONSOFMATRICES}.
\begin{lemma}
 Given $S, I \in \mathbb{C}^{n \times n}$ where $I$ is the identity matrix, $\sigma, \alpha \in \CC$ and $f$ is a complex analytic function, the following relation is fulfilled
 \begin{align*}
 f\left(
 \begin{bmatrix}
  \sigma I+\alpha S	&	I \\
  0	&	\sigma I
 \end{bmatrix}\right) =
  \begin{bmatrix}
  f(S)	&	f[\sigma I+\alpha S,\sigma I] \\
  0	&	f(\sigma)I
 \end{bmatrix}.
\end{align*}
\end{lemma}

\end{document}